\documentclass{article}
\def\tU{ \tilde{U} }

\def\bF{ \bar{F} }
\def\tI{ \tilde{I} }
\def\tJ{ \tilde{J} }
\def\phihat{ \hat{\phi} }
\def\TB{ \tilde{B} }
\def\Lie{ \mbox{Lie} }
\def\AF{ AF_{\bf Q} }
\def\tom{{ \tilde{\omega} }}
\def\d{ \delta }
\def\bA{ {\bf A} }
\def\bC{ {\bf C } }
\def\DGA{ DGA_{\bf Q} }
\def\pdr{ \pi^{dr}_1 }
\def\pcr{ \pi^{cr}_1 }
\def\A{{ {\cal A}_{\bf Q} }}
\def\B{{ {\cal B} }}

\def\Sh{{ Sh_{W,\bf Q}}}

\def\Spec{{ \mbox{Spec} }}

\def\Hom{{ \mbox{Hom} }}

\def\1ox{{ \Omega^1_{\scriptstyle{X}} }}
\def\2ox{{ \Omega^2_{\scriptstyle{X}} }}

\def\ok1{{ \Omega^1_K }}
\def\ok2{{ \Omega^2_K }}
\def\Om{{ \Omega }}
\def\om{{ \omega  }}
\def\O{{ {\mathcal O} }}

\def\ra{{ \rightarrow }}
\def\da{{ \downarrow }}

\def\a{{ \alpha }}

\def\e{{ \epsilon }}

\def\hra{{ \hookrightarrow }}
\def\da{{ \downarrow }}

\def\N{{ {\bf N} }}
\def\Q{{ {\bf Q} }}
\def\D{{ \Delta }}
\def\C{{ {\bf C} }}

\def\8{{ {\infty } }}

\def\s{{ {\sigma } }}
\def\G{{ \Gamma }}

\def\deg{{ \mbox{deg} }}
\def\^{{ ^{\wedge} }}

\def\Z{{ {\bf Z } }}

\newtheorem{thm}{Theorem}
\newtheorem{cor}{Corollary}
\newtheorem{lem}{Lemma}
\newtheorem{prop}{Proposition}

\newtheorem{definition}{Definition}
\usepackage{amsmath}
\usepackage{amsfonts}
\usepackage{amscd}
\usepackage{amssymb}

\def\invlim{\varprojlim}
\def\dirlim{\varinjlim}

\def\D{{ \Delta }}
\def\cD{ {\cal D} }

\def\hX{{ \hat{ X} }}

\def\cV{{ {\cal V} }}

\def\bcV{{ \bar{\cal V} }}
\def\D{{ \Delta }}
\def\tD{ \tilde{D} }

\def\H{{ \underline {H} }}
\def\C{{ {\cal C}_{W,\Q} }}

\title{A De\thinspace Rham-Witt 
approach to crystalline rational homotopy theory}
\author{ Minhyong Kim and Richard M. Hain}

\begin{document}
\maketitle
\begin{abstract} We give a definition of the
crystalline fundamental group of suitable log schemes
in positive characteristic using the techniques of
rational homotopy theory applied to the De\thinspace Rham-Witt
complex. 
\end{abstract}

\section{Introduction}
Throughout this paper, $k$ will denote a perfect field of
characteristic $p>0$, $W_n=W_n(k)$  the ring of
Witt vectors of $k$ of length $n$, and $W=\invlim W_n$.
$K$ denotes the fraction field of $W$.

We wish to present a relatively direct approach to the
unipotent crystalline fundamental group of a 
variety over $k$
using the De Rham-Witt (DRW) complex of Bloch and Illusie as complemented
by Hyodo and Kato (\cite{Bl}, \cite{Il}, \cite{Hyo1}, \cite{Hyo2},
\cite{HK})
and constructions that come from
rational homotopy theory (\cite{Su}, \cite{Q}, \cite{Ch}) 
and its Hodge-De Rham realizations (\cite{Mor}, \cite{Ha1}, \cite{Ha2},
\cite{N-A}). In the process,  to
a smooth connected proper fine log scheme $Y$ over $k$ of Cartier
type,
 we will associate a canonical commutative differential
graded algebra that deserves to be called the {\em unipotent
 crystalline
rational homotopy type}. The zero-th cohomology of the
bar construction
on this algebra will then give us the coordinate ring of the unipotent
 crystalline fundamental group. The use of the DRW complex allows
us to easily endow the crystalline fundamental group with natural
expected structures such as  Frobenius and monodromy operators and,
importantly, a weight filtration.

In the influential paper \cite{De1}, Deligne outlined
a motivic theory of the fundamental group. Given a variety
$V$ defined over a number field $F$, one should have
a  unipotent algebraic fundamental group corresponding to each cohomology
theory associated to $V$, Betti, \'{e}tale, de Rham, and
crystalline, together with suitable
comparison isomorphisms between them.  Following up on this
idea,
Deligne gave a Tannakian definition of a crystalline fundamental
group {\em using} the theory of the De Rham fundamental group: If $v$
is an absolutely unramified prime of $F$, and $F_v$ the completion  of $F$
at $v$, one assumed
that $V\otimes F_v$ had an integral model
$\cV$  with a smooth relative
normal crossing compactification $\bcV$.
This allowed him to define a Frobenius action
on the category of unipotent flat connections on
$V\otimes F_v$, and hence, on 
the unipotent de Rham fundamental group
$\pdr (V\otimes F_v)$ (which is isomorphic to  $\pdr (V) \otimes F_v$).
The de Rham fundamental group over $F_v$ together with this Frobenius
action was then taken to define the crystalline fundamental
group of the special fiber. This definition was reasonable in
the context of Deligne's paper since he was primarily
interested in varieties over number fields and,
of course, one could have started from
a variety defined over $F_v$ in this construction.
This  corresponds to the viewpoint that for a
variety over $F_v$, the crystalline cohomology is
just `extra structure' on De Rham cohomology.
However, it is evident that this definition is not quite satisfactory.
 The most important problem is that
 the crystalline fundamental group is not
 defined intrinsically for a variety in positive
characteristic. From the viewpoint of characteristic
$p$, it is  analogous to
defining the crystalline cohomology of a
variety  to be the de Rham
cohomology of a lifting (if it exists) to characteristic zero.
Another problem is to deal with the case of bad reduction.

In \cite{CL} and \cite{Sh1}, \cite{Sh2}, Chiarellotto and Le Stum,
and independently,
A. Shiho, gave an intrinsic definition of a crystalline
fundamental group $\pcr$ of a  proper log smooth
 variety over a perfect field
of positive characteristic $p$ which is
equipped with a comparison isomorphism
to $\pdr$ in a `lifted' setting. In Shiho, for example, the
the crystalline fundamental group was defined
as the fundamental group of the category of
unipotent isocrystals, and the  comparison theorem
was effected through the intermediary of a `convergent'
fundamental group, interesting in its own right.

In the present paper, we use an approach to the crytalline
fundamental group
suggested by
the work of Wojtkowiak \cite{Wo}  where the De Rham
rational homotopy groups were constructed
 over an arbitrary ground field using cosimplicial
schemes. Shiho's announcement \cite{Sh3} also indicates a similar approach
using a category of `complexes' of schemes.
The idea of both authors  is to  construct a model for
path space purely in the context of algebraic geometry
out of which one can extract the homotopy groups
via cohomological techniques.

On the other hand, we exploit the fact that
a formal algebraic model  already  exists for the cohomology
of path space, namely, the bar construction,
provided one has a commutative
differential graded algebra (henceforward, CDGA) which is a `purely algebraic'
analogue of the smooth differential forms. 
The problem here, as in rational
homotopy theory, is that the
usual resolutions which one uses to compute various
 cohomology
groups of sheaves of CGDA's
do not give rise to a CDGA in general, even though they
are  equipped with a multiplication which is homotopy
commutative.
However, it turns out that a CDGA  model with the right properties
 is available from the
work of Navarro-Aznar\cite{N-A} (already used in \cite{Wo}), where a derived
`Thom-Whitney' functor is constructed that associates in a canonical
fashion a CDGA to a {\em sheaf} of CDGA's on a topological space,
and more generally, a functor $R_{TW}f_*$ for
 maps $f:M \ra N$ of spaces, which converts
CDGA's on $M$ to CDGA's on $N$. This functor is equivalent to
the usual $Rf_*$ functor (defined on a suitable derived
category of sheaves) when we forget the multiplicative
structure.
Although the setting for Navarro's work is the Hodge
theory of complex varieties, it is clear that there is
actually a powerful technique that applies to quite
general topoi underlying his constructions.
More precisely, all one needs to apply his
machinery is the existence of acyclic Godement
resolutions of sheaves. In particular, one
obtains Navarro's functors for
CDGA's in the \'etale topology or the Zariski topology
of schemes. 
 
The
CDGA we use for our definition is the algebra of De Rham-Witt (DRW)
differential
forms, considered as a pro-sheaf on the small
\'etale site of a variety. This is very natural since
in many senses, the DRW differentials are the `correct'
analogue of algebraic differential forms suitable
for crystalline constructions.
This approach has the added advantage
that the definition of $\pcr$  is quite elementary and the
comparison with the De Rham fundamental group
requires only existing cohomological techniques.

But perhaps the main interest in the approach we present is
the natural definition of a weight  filtration on the
cohomology of the bar construction, and hence, on the coordinate
ring and the Lie algebra
of the crystalline fundamental group. For this, we use
a variant of the Hyodo-Steenbrink complex introduced by
Mokrane \cite{Mok}, modeled on the complex underlying
the `limit mixed Hodge structure' for homotopy groups \cite{Ha2}.

Referring to the following sections for precise terminology,
we state now the main theorems of the paper.

Let $Y$ be a connected proper smooth fine log scheme over  $k$ of Cartier type,
where $k$ is endowed with some fine log structure,
 and let $W\om_Y$ be the pro-sheaf of CDGA's consisting
of the De Rham-Witt differential forms of Illusie-Hyodo-Kato.

The {\em crystalline rational homotopy type}
 $A_Y$ of $Y$ is defined by the following
formula:
$$A_Y=s_{TW}\invlim_K \G (G (W\om_Y))$$
The notation, which will be explained below in detail,
is that given a pro-sheaf $L$,
$G(L)$ is its canonical cosimplicial Godement resolution,
$\invlim_K$ is the operation which associates to an inverse
system of (cosimplicial) $W_n$-modules the (cosimplicial) 
$W$-module obtained by taking
the inverse limit, and then tensors it with $K$, and finally,
$s_{TW}$ is Navarro-Aznar's `simple Thom-Whitney algebra'
functor.
A choice of a point $y$ determines an augmentation for
$A_Y$, and we can form the bar complex, which we denote
by $B(Y,y)$. Then $Cr(Y,y):=H^0(B(Y,y))$ has the natural structure of
a commutative Hopf algebra filtered by finitely
generated Hopf algebras.
The {\em crystalline fundamental group} is defined by
 $$\pcr (Y,y)=\Spec Cr(Y,y)$$
When $k$ is equipped with the
log structure of a `punctured point'
and the point gives an exact embedding of log schemes, 
we will see below that $\pcr(Y,y)$ is naturally equipped
with a semi-linear Frobenius as well as a monodromy operator,
which is interpreted as a vector field on $\pcr$, and that
they
satisfy the usual relations.

When $Y$ is semi-stable, the weight filtration
on the DRW complex induces one on the bar complex and
on $Cr(Y,y)$.
An obvious formalism of `mixed Frobenius complexes'
then yields our main result:
\begin{thm}
Assume $k$ is  finite and that $Y$ is globally the union of
smooth components that meet transversally.
Then the spectral sequence for the weight filtration
degenerates at $E_2$ and gives $Cr(Y,y)$ the
structure  of a mixed isocrystal.
\end{thm}

This structure is compatible with the Hopf algebra structure,
and hence, also induces the structure of a mixed isocrystal on
$\Lie (\pcr(Y,y))$.

If $X$ is a log smooth variety over a field $F$ of characteristic 0,
we  associate to $X$ its De Rham rational homotopy type
 by the formula
$$A_X=s_{TW}\G (G (\Om_{X/F}))$$
Given an augmentation associated to a point $x$,
one then forms the bar complex $B(X,x)$
and defines the De Rham algebra $$DR(X,x):=H^0(B(X,x))$$
The De Rham fundamental group is
given by \cite{Wo}
$$\pdr(X,x)=\Spec (DR(X,x))$$
As mentioned previously, one advantage of our approach is
that Berthelot-Ogus-type comparison theorems are within
the scope of crystalline cohomological techniques.

Let $A$ be a complete local ring
with fraction field $F$ and perfect residue field $k$.
Let $W$ be the ring of Witt vectors
of $k$ with fraction field $K$.

\begin{thm}
Suppose $X$ is a proper connected smooth  scheme over $A$
 with a relative
normal crossing divisor $D$.
Equip $X$ with the log structure associated to the
divisor $D$.
 Denote by $X^*$ the generic fiber of $X$ (with the induced
log structure),
$x$ a point of $X^*-D $ with reduction $y\in Y-D$. Then
$$A_{X^*}\simeq A_Y\otimes_K F$$
where the isomorphism is in the
homotopy category of commutative differential graded
algebras over $F$. Furthermore,
$$Cr(Y,y)\otimes_K F \simeq DR(X^*,x)$$
as commutative Hopf algebras over $F$.
\end{thm}
\begin{cor}With the assumptions of the theorem,
$$\pdr(X^*,x) \simeq \pcr(Y,y)\otimes_K F$$
\end{cor}

Although the main emphasis in this paper is on fundamental groups,
the crystalline rational homotopy type can be used
to define higher crystalline homotopy groups in the
simply-connected case. In particular,
we also get some interesting consequences of the Artin-Mazur
type for the higher rational homotopy groups of simply connected
varieties over number fields from this comparison theorem.

\begin{cor}
Let $X$ and $X'$ be smooth proper
  connected varieties over a number field $F$ 
equipped with normal crossing divisors $D$ and $D'$.
Let $v$ be a prime of $F$ where both $(X,D)$ and
$(X',D')$ have  good reduction in the sense that
both varieties extend to smooth proper
schemes over $\O_{F_v}$ and the divisors extend to
relative normal crossing divisors.
Suppose
the pairs are isomorphic mod $v$
and both $X-D$, and $X'-D'$ are simply-connected
in any embedding of $F$ into the
complex numbers.
Then the higher rational homotopy groups 
 of $X-D$ and $X'-D'$
are isomorphic for any embedding of
$F$ into the complex numbers.
\end{cor}

Wojtkowiak \cite{Wo} had earlier shown that these groups
are also independent of the embedding.

One issue that is completely ignored in this paper
is the comparison with the Tannakian view. That is,
we do not show that the crystalline $\pi_1$
we define classifies unipotent isocrystals.
We hope to carry out this comparison  in a
subsequent paper. Another topic
we hope to deal with is relative
completions of crystalline fundamental groups with
coefficients in an $F$-isocrystal,
 the crystalline analogue of the completion of the
fundamental group of a smooth variety relative to a variation of Hodge
structure, which is considered in \cite{Ha3}.
For a related but different approach, the reader is also referred
to the paper of Vologodsky \cite{Vo}.

In a forthcoming publication, we will present
a generalization of the comparison
isomorphism to incorporate the `Hyodo-Kato case' of
semi-table reduction over a ramified base.
This will be achieved by using an `infinitely-twisted
telescope' construction and the ideas of
\cite{Og}. These ideas will also be applied to
prove a $p$-adic analogue of Oda's
good reduction criterion for curves
\cite{Od}.

\section{Review of Hodge-De Rham theory for homotopy groups}
The unipotent De Rham fundamental group  
$\pdr (X,x)$ of a space $X$, say  with coefficients in $\bC$,
can be defined as the complex pro-unipotent completion of
the usual fundamental group $\pi_1(X,x)$. That is, $\pdr(X)$
is the initial object in the category of
 inverse systems of pro-algebraic unipotent
groups $U$ over $\bC$ equipped  with group homomorphisms $\pi_1(X,x) \ra U$.

Assuming the space has
finite-dimensional  $H^1$,
one realization is constructed  by considering the group
algebra $R=\bC \pi_1(X,x)$ together with the augmentation ideal $J$,
and then the completion $\hat{R}$ of 
$R$ w.r.t. the augmentation ideal. $R$ naturally has the structure
of a Hopf algebra induced by the comultiplication $c$ defined on the image
of elements of $\pi_1$ by
$c(g)= g \otimes g$. This comultiplication extends uniquely to
 $\hat{R}$ and the complex points of
$U$ can be realized as the
group-like elements in $\hat{R}$, i.e., $u \in \hat{R}$
such that $c(u)=u\otimes u$.

More precisely, $U$ is defined by  the inverse system
 given by group-like elements in $\hat{R}/\hat{J}^n$.
The Lie algebra of the De Rham fundamental group
can then be realized as the primitive elements inside $\hat{R}$,
i.e., those elements $t$ that satisfy
$c(t)=t\otimes 1 + 1 \otimes t$.

Another way of understanding this construction
is to consider the dual ind-Hopf algebra 
$$R^*:=\dirlim \Hom_{{\bf C}}(R/J^n, {\bf C}).$$
The remarks above correspond to the fact that
 $R^*$ is the affine coordinate ring of
$\pdr (X,x)$

If $X$ is a smooth complex variety, $\pdr(X,x)$ is also the Tannaka
group of the category of unipotent vector bundles with
flat connection associated to the fiber functor of evaluation at
$x$.

The point of view that informs this paper comes
from the construction of the complex (or real) De Rham fundamental group
 via iterated integrals of
 differential forms. That is, if $X$ is a manifold,
consider the CDGA $A_X$ of $C^{\infty}$ differential forms on $X$
with the augmentation $a_x$ given by pull-back to the point $x$.

In Chen's approach \cite{Ch}, one constructs locally constant
functions on the
loop space at $x$ via iterated integrals, which suffices to
construct the coordinate ring. This process is conveniently
formalized using the bar complex $B(A_X, a_x)$
associated to the augmented algebra. The detailed definition will
be reviewed below, but we note that the essential
part of the zero-th degree term
is given by 
$$\oplus^{\infty}_{s=0} (\otimes^s A^1_X)$$
Given a tensor product $a_1\otimes \cdots \otimes a_s $ from this
vector space, we get a function
 on loop space according to the rule
$$\gamma \mapsto \int_{\gamma} a_1 a_2 \cdots a_s$$
where the last quantity is the iterated integral
defined according to the following prescription:
Write $\gamma^*(a_i)=f_i(t)dt$. Then
$$\int_{\gamma} a_1 a_2 \cdots a_s:=\int_{0\leq t_1\leq t_2 \leq \cdots \leq t_s \leq 1}
f_1(t_1) f_2(t_2) \cdots f_s(t_s) dt_1 dt_2 \cdots dt_s$$
The elements of $H^0(B)$ then correspond to locally
constant functions on the loop space. 

We will see below that $B(A_X, a_x)$
also
has the structure of a Hopf algebra which induces a commutative
Hopf
algebra structure on $H^0(B)$.

Chen's theorem says that
$$\pdr (X,x)_{{\bf C}}\simeq \Spec (H^0(B))$$
provided we interprete the right hand side suitably
as a pro-algebraic group. (This is done using the
bar filtration, also reviewed below.)
In the simply-connected case, we can also recover the
higher complex homotopy groups from the higher cohomology
of the bar complex.

One advantage of this approach as opposed to Sullivan's
theory of minimal models is that the transparent relation
between $A_X$ and the bar complex enables one to
carry over extra structure from the differential
forms naturally to the bar complex. In precise terms,
the bar complex is {\em functorial} in $A_X$, while the
minimal model is functorial only in an apropriate
homotopy category.
This allowed, for example,
the construction of natural mixed Hodge structures on the
coordinate ring of $\pdr (X,x)$ \cite{Ha1}  as well as on
$\Lie (\pdr(X,x))$ for general varieties over ${\bf C}$.
Equally important is the fact that the bar construction is
completely algebraic, once one is given $A_X$. Thus, it
can be built on an arbitrary CDGA over any field of characteristic zero, 
in
contrast to the iterated integrals which motivate its definition
but are specific to differential forms on a manifold.
It is this algebraicity which allows us to  to give a natural
construction of the crystalline fundamental group together
with all the extra structures it should be endowed with.
We will carry  this out in subsequent sections.

\section{Algebraic prerequisites: The Thom-Whitney functor and the bar complex}
We will quickly review the definitions of \cite{N-A} and \cite{Ha1}. 

Given a category $T$, we denote by $\D T$ 
the category of cosimplicial objects
in $T$ and by $\D^+ T$ the
category of augmented cosimplicial objects in
$T$ \cite{De3}.

Let $Y$ be a scheme over $k$.

 Denote by $\A (Y)$ the
category of  pro-objects
 in the category of sheaves of
(graded-)commutative differential-graded $W$-algebras (CDGA's)
`up to isogeny'
on the small \'etale site of $Y$. 
So an object of $\A (Y)$ consists of a sequence
$A=(A_n)_{n\geq 1}$, where each $A_n$ is a $W_n$-algebra
and we are given transition maps $A_{n+1}\ra A_n$ which
we assume to be a surjection.
 Morphisms from $A=(A_n)$ to $B=(B_n)$
are given by
$$\Hom (A,B)=[\invlim_j \dirlim_i \Hom (A_i, B_j)]\otimes_{\Z} \Q$$
We do not distinguish notationally between the ring $W$
and the pro-sheaf $W$ regarded as an object of $\A (Y)$ for various $Y$.
Also, the scheme $Y$ will often be suppressed  from the
notation if the context makes it unnecessary.

Note that an object $A$ of $\A$ also carries a superscript
corresponding to the complex degree and an object
$G$ of $\D \A$ carries two superscripts $G^{\bullet ,*}$, the
first denoting cosimplicial degree and the second complex
degree. These scripts will usually be suppressed from our
notation. In other contexts as well,
we will avoid making scripts explicit unless necessary for
clarity. (In fact, it is our experience that it enhances
clarity most of the time to avoid thinking about scripts.)

A word about our convention:
When referring to various (pro-)sheaves with extra structure,
we will often suppress the `sheaf' in their designations unless
forced upon us by considerations of clarity, and
similarly be somewhat careless
with localization up to isogeny. So an object
of $\A (Y)$ will usually just be called a CDGA on $Y$.
Also, unless explicitly stated otherwise, an algebra will refer
to a $W$-algebra.
On the other hand, 
when we put $k$ into the argument of one of our categories, such as
$\A(k)$, we will be referring to actual objects and not sheaves.
Thus, in the previous sentence, we mean the category of
inverse systems of CDGA's over $W$, and not the category
of sheaves on the \'etale site of $\Spec(k)$.
Various constructions will be described in
the sheaf case and can be modified in an obvious way for
objects over $k$.

In a manner analogous to $\A (Y)$, we define $\Sh (Y)$, the
category of pro-sheaves of $W$-modules up to isogeny,
and $\C (Y)$, the category of pro-complexes of sheaves of $W$-modules up to
isogeny, with similar conventions of reference as explained
in the previous paragraph. Finally, $\DGA (Y)$
will be the category of pro-differential graded $W$ algebras
up to isogeny. Thus, $\A (Y)$ is a full subcategory of
$\DGA(Y)$ which, in turn, is equipped with a natural forgetful functor
to $\C (Y)$.

Given an object $C$ of $\C$, its cohomology sheaves are by
definition $\H^i(C):=(\H^i(C_n))$.
In what follows, of particular importance will be the situation
where the $\H^i(C) $ are objects of
$\Sh$, that is, where the transition maps are
surjective. Sometimes, we will denote by $\H(C)$ the direct sum of the
$\H^i(C)$ considered as a complex with zero differential.
Thus, if $A\in \A$ and the transition maps for cohomology are
surjective, then $\H(A)$ is naturally an object of
$\A$ with zero differential. We will circumvent some foundational
annoyances by passing to cohomology sheaves only in situations
where they belong to $\Sh$.
In that case, given two objects $C$ and $C'$ of $\C$, we say they are
quasi-isomorphic (q.i.) if there is a map 
$f:C\ra C'$ in $\C$
which induces an isomorphism (in $\Sh$) of cohomology pro-sheaves.
We say $A$ and $A'$ of $\A$ are multiplicatively q.i. if
they are q.i. as objects of $\C$ but via a map in the category
$\A$.

The homotopy category of  $\A$ is obtained by keeping
the same objects and inverting all the multiplicative quasi-isomorphisms. 
We will
say two objects in $\A$ are quasi-equivalent (q.e.) if they
are isomorphic in the homotopy category. 

Denote by $\bC (K)$ the category of complexes of $K$ vector
spaces and $\bA (K)$ the  category of CDGA's over $K$,
equipped with
the natural forgetful functor to $\bA(K)$.

The functor $\invlim_K$ from $\C (k)$ to $\bC(K)$ takes
a pro-complex of $W_n$-modules to its inverse limit,
which is a $W$-module, and then tensors with $K$.
$\invlim_K$ takes $\A(k)$ to $\bA (K)$. We use the
same notation for the functor induced on $\D \C(k)$.

As explained in \cite{N-A}, pp. 13-14 and p. 23, given
an object $G$ of $\D \bA (K)$ there is the `usual' way
of  giving to $s(G)$, the associated simple complex,  
a multiplication, making it into a DGA
over $K$, depending on the
 choice of  `Eilenberg-Zilber transformations.'
We will use the notation $s(G)$ for this DGA, whether or not
we are remembering its multiplicative structure. 
This multiplication will not be graded-commutative in general.
On the other hand,
the multiplication induced on
$H(s(G))$ is commutative, and it can be considered
as an object of $\bA(K)$.

To construct the
Thom-Whitney algebra requires a choice of
an algebraic model for differential forms on the
standard simplices. 
Let $E_n$ be the algebra of global differential forms relative to $\Spec (K)$
on the variety $\D_n:=\Spec (K [t_0, \ldots, t_n]/(\sum t_i-1))$.
The $\D_n$'s form a cosimplicial scheme in the usual manner via
co-face maps $\d^i:\D_n  \ra \D_{n+1}$ 
given by
$$(t_0, \ldots, t_{n}) \mapsto (t_0, \ldots, t_{i-1},0, t_{i},
\ldots t_{n})$$
and codegeneracy maps
$s_i: \D_{n+1}\ra \D_n$:
$$ (t_0, \ldots, t_{n+1}) \mapsto (t_0, \ldots, t_{i}+t_{i+1},
\ldots t_{n+1})$$
 so that the $E_n$'s form a simplicial
CDGA over $K$ that we will denote by $E$. A standard
computation shows that $s(E)$ only has cohomology in
degree 0 of dimension  1.

Now, the Thom-Whitney algebra $s_{TW}(G)$ of $G$ is by definition
the simple CDGA  associated to the
end \cite{Mac} of $E\otimes_{K} G$ considered as a functor from
$\D^{op}_{mon} \times \D_{mon}$ to  commutative differential bi-graded
algebras over $K$. Here, $\D_{mon}$ refers to the
subcategory of $\D$ where the morphisms are
 strictly increasing maps.
Thus, it is an object of 
$\bA (K)$. We elaborate a bit on this definition (using scripts):
Elements of $E\otimes G$ will be of the form 
$$(\sum e^p_n\otimes g^{m,q})_{n,m}$$
where $n$ is the simplicial degree, $m$, the cosimplicial
degree, and $p,q$ are the complex degrees. The elements of the end 
are compatible sequences (indexed by $n \in \N$)
$$[\sum e^p_n\otimes g^{n,q}]_{n}$$
where compatibility refers to the equality:
$$ \partial_i\otimes 1 [\sum  e^p_{n+1}\otimes g^{n+1,q}]=1\otimes \d^i
 [\sum e^p_n\otimes  g^{n,q}] \in E_n \otimes G^{n+1} $$
for all $n$. Since we still have the complex degrees $p$ and $q$
left, the result is a CD bi-graded algebra. One then takes
the associated simple complex to get the  CDGA $s_{TW}(G)$.

Readers unfamiliar with the notion of ends should apply
it to the bifunctor $E_.\otimes S^.$ as an exercise,
where $S^.$ is the cosimplicial algebra (not the complex)
of singular cochains with $K-$values on a topological space
$X$. (One gets Sullivan's polynomial differential forms on $X$
with values in $K$.)

We note that both
$s$ and $s_{TW}$  define functors $s,s_{TW}:\D \bA (K)  \ra \bC (K) $.
By using an integral version of the simplices, we could
define $s_{TW}$  also for cosimplicial $W$-algebras so that
the construction commutes with $\invlim_K$. On the the other hand,
the comparison with $s$, to be discussed below, will
not work at the integral level. 

Given choices of geometric points over all the points of
$Y$, we can construct cosimplicial Godement resolutions in the usual
fashion level by level \cite{Go} for a prosheaf.
So from an object $C\in \C (Y) $ we obtain  an object of $\D \C (Y)$, which we will
call the Godement resolution of $C$ and denote by $G(C)$.
If $A\in \A $, then $G(A) \in \D \A$.
Applying the global section functor $\Gamma$ and the associated
simple complex functor $s$, we get an object
$s (\G (G(C))) $ in $\C(k)$ which represents
$R\G (C)$. When we write $R\G (C)$, therefore, we
will  mean this explicit pro-complex.
More important for our purposes is
$$R\G_K(C):=\invlim_K s (\G (G(C)))=s(\invlim_K (\G (G(C))))$$
It is a standard fact that if $C$ is quasi-isomorphic to $C'$,
then $R\G_K(C)$ is q.i. to $R\G_K(C')$.

As mentioned above, if we start with
$A \in \A (Y)$,  $R\G_K (A)$ will only be a DGA over $K$,
not necessarily commutative. On the other hand,
we can construct the object
$$TW(A):=s_{TW}(\invlim_K \G ( G(A))).$$ 
which lies in $\bA (K)$.

The proof of the following  are in Navarro-Aznar's paper.

\begin{lem}[\cite{N-A} Thm. 2.14]
There is a natural transformation of functors on
$\D \A (k)$,
$$ I: s_{TW}\circ \invlim_K \ra s \circ \invlim_K  $$
which induces
quasi-isomorphisms when evaluated on objects of
$ \D \A (k)$, and furthermore, induces an isomorphism
of CDGA's at the level of cohomology.
\end{lem}

\begin{cor}[\cite{N-A}, (4.4)] 
Given $A\in \A (Y)$, there 
is a q.i. of complexes 
$$TW(A) \simeq  R\G_K(A)$$
which induces an algebra isomorphism on cohomology.
\end{cor}

\begin{cor}[\cite{N-A} (4.6)]
If $A$ is q.i. to $A'$ in $\A (Y)$, then
$TW(A)$ is q.i. to $TW(A')$ as CDGA's over $K$.
\end{cor}

\begin{lem}[\cite{N-A} (3.7)]
Let $C^+$ be an object of $\D^+ \A (k)$ with $C^{+,-1}=A$ and $C:=C^+|\D$.
Thus we have two maps of complexes
$$\invlim_K A\ra s_{TW}\invlim_K(C)\ \ \ \ \invlim_K A \ra s \invlim_K(C)$$
The natural tranformation $I$, gives rise to a commutative
diagram:
$$\begin{array}{ccc}
& & s_{TW}\invlim_K(C) \\
& \nearrow & \da \\
\invlim_K A & \ra & s \invlim_K(C)
\end{array}$$
Thus, if the second map is a q.i., then the first is a multiplicative
quasi-isomorphism.
\end{lem}

It will be convenient to have the Thom-Whitney functors
also defined for co-simplicial sheaves of
CDGA's. This is
easily achieved by applying $s_{TW}$ twice:
If $A^{.}$ is a cosimplicial CDGA on $Y$,
then each 
$TW(A^n)$ is a CDGA over $K$ and they fit together
to form a cosimplicial CDGA over $K$. Applying
$s_{TW}$ to this gives us a CDGA
that we will denote by
$TW(A^{.})$.
The following is easily deduced by integrating twice:
\begin{lem}
$$TW(A^{.}) \simeq R\Gamma_K (s(A^{.}))$$
\end{lem}

Similarly, suppose $Y_.$ is a simplicial scheme and $A_.$ is
a CDGA over $Y$. We see then that the $TW(A_n)$ are
objects of $\bA (K)$ and they come together to
form an object of $\D \bA (K)$. We simply apply
$s_{TW}$ again to get
$$TW(A_.):=s_{TW}(\{TW(A_n)\}_n)$$
The usual $R\G_K$ on such an object can be constructed
as $s(\{R\G_K (A_n)\}_n)$
so applying the integration functor twice gives us
a q.i.
\begin{lem}
$$TW(A_.)\simeq R\G_K(A_.)$$
\end{lem}

\begin{cor}
Let $p:Y_. \ra Y$ be a simplicial  hypercovering which satisfies
cohomological 1-descent
for the \'etale topology. Let $A$ be a CDGA on $Y$, 
$A_.$ a CDGA on $Y_.$, and suppose we have a q.i.
$p^*A \simeq A_.$.
Then this induces a q.i.
$$TW(A) \simeq TW(A_.)$$
\end{cor}

There is a version of the functor $TW$
for filtered CDGA's (\cite{N-A} section 6).

Denote by $\AF (Y)$ the category of filtered pro-CDGA's on $Y$
up to isogeny. So an object is a pair $(A,F)$ where
$A$ is an object of $\A (Y)$ and $F$ is a multiplicative
decreasing
filtration of $A$.
Thus, we are given subobjects
$F_n^iA_n$ for each $i\in \Z$ and each level $n$
such that $F_n^{i+1}A_n \subset F_n^iA_n $ and the
transition maps send $F_n^iA_n$ surjectively to $F_{n-1}^iA_{n-1}$.
Furthermore, we have $F_n^iF_n^j\subset F_n^{i+j}$ at each level.
 Finally, the morphisms are defined by
$$\Hom_{\AF (Y)}((A,F),(A',F')):=[\invlim_m \dirlim_n 
\Hom ((A_n, F_n), (A'_m,F'_m))]\otimes \Q$$
Let $(A,F)$ be an object of $\D \AF (Y) $. Then we can define the
filtered CGDA
$TW (A,F)$ over $K$.
in a manner entirely analogous to the previous discussion: One takes the simple
filtered CDGA associated to the end of
$E\otimes \invlim_K \G (G(A,F))$ where the filtration on the tensor product is
induced by the given filtration $F$ on $A$ and the trivial decreasing
filtration $\e$ of $E$ defined by 
$$\e^i(E)=
\left\{ \begin{array}{ccc} E & , & i\leq 0 \\ 0& , & i>0 \end{array}
\right.$$

We also have the functor $s$ which associates to $(A,F)$the filtered
DGA $(s(K), s(F))$.
When applied to objects of $\D \AF (k)$, we can again compose
with inverse limits to end up with cosimplicial filtered CDGA's
over $K$.

\begin{lem}[\cite{N-A} Lemma 6.3, (6.7)]
The natural transformation  $$I:  s_{TW}\circ \invlim_K \ra s\circ \invlim_K $$
of functors on $\D \AF (k) $
discussed previously  also gives a filtered quasi-isomorphism which induces
in cohomology an isomorphism of filtered CDGA's.
\end{lem}

If we denote by $R\G_K (A,F)$ the filtered DGA $$(R\G_K (A), s(\invlim_K\G (G(F)))) $$ 
Then
\begin{cor} There is a filtered q.i.
$$TW(A,F) \simeq \invlim R \G_K (A,F)$$
of complexes of $K$-vector spaces which induces an isomorphism
of filtered CDGA's in cohomology.
\end{cor}

\begin{cor}
$TW$ takes filtered q.i.'s to filtered q.i.'s (\cite{N-A} (6.14))
\end{cor}

We now give a brief discussion of the bar complex.

Let $A\in \bA(K)$ and let $a:A \ra K$ be an augmentation. Assume
throughout this discussion that $A$ has connected cohomology,
that is, $H^0(A)=K$, and has no cohomology in negative
degrees.
Define the  bar complex $B(A,a)\in \bC (K)$ by 
the formulas in \cite{Ha1} pp. 275-276:
Let $I$ be the augmentation ideal
and
$$B^{-s,t}(A,a)=(\otimes^s I)^t$$
where the outer superscript $t$ denotes the subset
of elements of degree $t$. We denote the
element
$$ a_1 \otimes \cdots \otimes a_s$$
of $B^{-s,t}(A,a)$ by $[a_1|\cdots |a_s]$.

There are two differentials $d_C$ and $d_I$ (the `combinatorial'
and `internal' differentials)
$$d_C: B^{-s,t} \ra B^{-s+1,t}\ \ d_I:  B^{-s,t} \ra B^{-s,t+1}$$
given by the formulas:
$$\begin{array}{rl}
d_C([a_1|\cdots |a_s]):=&  \sum_{i=1}^{s-1}(-1)^{i+1}[Ja_1| \cdots |Ja_{i-1}|Ja_i \wedge a_{i+1}| a_{i+2}\cdots |a_s] \end{array}$$
 where $J(v)=(-1)^{\deg v}v$
and
$$\begin{array}{rl}
d_I([a_1|\cdots |a_s]):=& \sum_{i=1}^{s}(-1)^i[Ja_1| \cdots |Ja_{i-1}|da_i|a_{i+1}| a_{i+2}\cdots |a_s]\end{array}$$

These differentials make
the direct sum of the $B^{-s,t}$ into a double complex
and we denote by $B(A,a)$ the associated total complex.

There is a filtration $\B$ on $B$ given by
$$\B^{-s}:=\oplus_{u\leq s} B^{-u,t}$$
 called the {\em  bar filtration}.

 Denote by 
$$\B_0 H(B(A,a))\subset \B_1 H(B(A,a)) \subset \B_2 H(B(A,a)) \subset \cdots$$
the filtration induced on the cohomology, where one notes
that the signs of the indices
 have been reversed with respect to the bar filtration.
That is,
$$\B_rH:=\mbox{Im}(H(\B^{-r}))$$

Let $E_r$ be the spectral sequence, the {\em Eilenberg-Moore} spectral
sequence, associated with the bar filtration and converging
to the cohomology of the bar complex.
We have for the $E_1$-term,
$$E_1^{-s,t}\simeq B^{-s,t}(H(A))$$
the terms of the bar complex on the cohomology of $A$ regarded
as an augmented algebra.

In particular,
$$E_1^{-s,s}= \otimes^s H^1(A)$$

In any case, we see that if $A $ has  coherent cohomology,
that is, if all the $H^i$ are finite-dimensional $K$-vector spaces,
then the  $\B_rH^n(B(A,a))$ are also finite dimensional for
all $n$. We note that for coherence of
 $\B_rH^0(B(A,a))$ we just need the coherence of
$H^1(A)$.
If $A$ has coherent cohomology which furthermore
occurs only in bounded degrees,
we also see that $\B_rH(B(A,a))$ is finite-dimensional
over $K$ for each $r$.

There is a multiplicative structure on $B(A,a)$ given by
the formula
\cite{Ha1} p.278:
$$[a_1|\cdots |a_r][a_{r+1}|\cdots |a_{r+s}]=\sum \e (\sigma)
[a_{\sigma(1)}|\cdots |a_{\sigma{r+s}}]$$
where $\sigma$ runs over all shuffles of type $(r,s)$and $\e:
\Sigma_{r+s}\ra \{\pm 1 \}$ is the representation
of the symmetric group obtained by giving
$a_j$ weight $-1+\mbox{deg}a_j$. That is, transpositions contribute
a minus sign only when switching two elements of even degree. 
 
We also have the comultiplication:
$$[a_1|\cdots |a_s]\mapsto \sum_{i=1}^s [a_1|\cdots|a_i]\otimes [a_{i+1}|\cdots|a_s]$$
that combines with the multiplication to give
$B(A,a)$ the structure of a differential graded Hopf algebra
over $K$. This induces a commutative Hopf algebra structure
on $H^0(B(A,a))$.
The bar filtration on $B(A,a)$ is preserved by
the comultiplication and hence, we get a filtration
of $B(A,a)$ by sub-Hopf-algebras $B_r(A,a)$, defined
to be the subalgebra generated by the $r$-th level of
the bar filtration. We will be mostly interested in
the induced filtration   of $H (B(A,a))$ which we will
denote by $H_r(B(A,a))$.
If $H (A)$ is coherent with bounded degree and $H^0(A)=K$, we see that this
  induces a filtration of $H (B(A,a))$
 by finitely generated subalgebras. Similarly, if
$H^1$ is finite-dimensional, then $H^0(B(A,a))$ is
filtered by finitely generated subalgebras $H^0_r(B(A,a))$.

We will need an algebraic fact about the bar complex
to compare fundamental groups to homology. The following
result is of course well-known (cf. e.g. \cite{Su}), but
we were unable to locate a proof involving the
maps that we need for our purposes. Therefore, we include one here.

As mentioned, with our assumptions,
$H^0(B(A,a))$ has the structure of a
non-negatively graded, commutative Hopf algebra.
 Each closed element
$f \in A^1$ determines a closed element $[f]$ of $H^0(B(A,a))$.
Since $[df] = d[f]$ for all $f$ in the
augmentation ideal,  $[f]$
depends only on the class of $f$ in $H^1(A)$. There is therefore a
well defined linear mapping
$$
\phi : H^1(A)[1] \to H^0(B(A,a)).
$$
Denote the bicommutative Hopf algebra generated by the vector space $V$
by $SV$. Since $H^0(B(A,a))$ is commutative, $\phi$ induces an
algebra homomorphism
$$
\phihat : S(H^1(A)[1]) \to H^0(B(A,a))
$$
Since each $[f]$ is a primitive element of $H^0(B(A,a)$, $\phihat$ is
a Hopf algebra homomorphism.

\begin{prop}
The homomorphism
$$
\phihat : S(H^1(A)[1]) \to H^0(B(A,a)).
$$
is the inclusion of the unique maximal cocommutative Hopf subalgebra.
\end{prop}

{\em Proof.}
First note that $S(H^1(A)[1])$ is naturally graded by the symmetric
powers of its primitives:
$$
S(H^1(A[1]) = \oplus_{s \ge 0} S^s (H^1(A)[1]).
$$
If we set
$$
B_s S(H^1(A)[1]) = \oplus_{t \le s} S^t (H^1(A)[1])
$$
then $S(H^1(A)[1])$ is isomorphic to its associated graded and
also the homomorphism $\phihat$ is filtration preserving. One therefore
has a homomorphism
$$
Gr^B \phihat : S(H^1(A)[1]) \to Gr^B H^0(B(A,a))
= \oplus_{s\ge 0} E_\infty^{-s,s}
$$

When $1 \le r < r' < \infty$, we have inclusions.
$$
\oplus_{s\ge 0} E_r^{-s,s} \supseteq \oplus_{s\ge 0} E_{r'}^{-s,s}
\supseteq \oplus_{s\ge 0} E_\infty^{-s,s}.
$$
Define $\phihat_r$ to be the composite
$$
S(H^1(A)[1]) \stackrel{Gr\phihat}{\longrightarrow}
\oplus_{s\ge 0} E_\infty^{-s,s} \hookrightarrow
\oplus_{s\ge 0} E_r^{-s,s}.
$$

Since each $E_r$ is a connected, graded Hopf algebra, and since $\phihat_r$
is injective on primitives, each $\phihat_r$, $1\le r \le \infty$, is
injective. Since $Gr^B\phihat$ is injective, $\phihat$ is also injective.

It remains to prove maximality. First, the dual of
$$
\oplus E_1^{-s,s}
$$
is the tensor algebra on the dual of $H^1(A)[1]$; its Hopf algebra
structure is characterized by the fact that the dual of $H^1(A)[1]$ is
primitive. This has, as maximal commutative quotient, the free bicommutative
Hopf algebra generated by the dual of $H^1(A)[1]$. It follows that
$\phihat_1$ is the inclusion of the unique maximal bicommutative sub Hopf
algebra of $\oplus E_1^{-s,s}$. But since
$$
\bigoplus E_1^{-s,s} \supseteq \oplus E_2^{-s,s} \supseteq \bigoplus
E_3^{-s,s} \supseteq \cdots \supset E_\infty^{-s,s}
$$
it follows that $\phihat_r$ is the inclusion of the unique maximal
bicommutative Hopf subalgebra of $\oplus E_r^{-s,s}$ whenever
$1 \le r \le \infty$.

The result now follows as $\phihat_\infty = Gr^B \phihat$ and since
$$
H^0(B(A,a)) = \bigcup_{s\ge 0} B_s H^0(B(A,a)).
$$

\section{Review of the De Rham-Witt complex}

We will work in the setting of log schemes.
and make extensive use of the theory developed
by Kato \cite{Ka}, Hyodo \cite{Hyo1} \cite{Hyo2},
and Hyodo-Kato \cite{HK}. The reader
should consult these articles for precise notions and notation.
Furthermore, a number of important gaps in the
literature have been filled in the recent
preprint of Nakkajima \cite{Nak}.
We will cite the necessary results as we proceed.

Denote by $S_0$ the scheme $\Spec(k)$
endowed with a fine log structure $L$. Unless the context
makes it necessary to be careful with the distinction,
we will denote by the same letter the scheme without
the log structure.
 $L$ determines
a canonical log structure on $W=W(k)$ 
 induced by the pre-log structure that composes
 $L\ra k$ with the Teichm\"uller lift.
One has a similar construction for any log scheme
$Y$ over $k$. By $W(Y)$, we denote the system of log schemes
(or ind-log scheme)
with underlying space the same as the space of $Y$,
but with structure sheaves $W_n(\O_Y)$ and the log structures
lifted with the Teichmueller character. This
definition extends naturally to simplicial log schemes
$Y_.$ over $k$ to give simplicial ind-log schemes
$W(Y_.)$.

Denote by $S$ the scheme $\Spec(W)$,
again with and without the log structure.
Let $Y$ be a  smooth, fine log scheme over $S_0$ of
Cartier type. For each $n\geq 1$, Hyodo and Kato define a
level $n$ De Rham-Witt complex $W_n\om_Y $ (w.r.t. $S_n$)
which is  a (sheaf of)  commutative differential
graded algebra(s)  on the small \'etale site of (of the underlying scheme of)
$Y$  equipped with
projections $\pi_n:W_n\om \ra W_{n-1}\om $. The log
de Rham complex $\Om_{Y/S_0}$ occurs 
at the bottom level (that is, $W_1\om_Y $) and in
degree zero, we have 
 $W_n\om_Y^0=W_n\O_Y$, Serre's sheaf of Witt vectors.
There are also operators
$$F:W_{n+1} \om^q \ra W_n \om^q, \ \ \ V: W_{n} \om^q \ra W_{n+1} \om^q $$
extending the usual Frobenius and Verschiebung on $W\O$
and satisfying the identities 
$$F(ab)=F(a) F(b), \ \ V(F(a)b)=a V(b)$$
and
$$FV=VF=p  \ \ \  FdV=d$$
which implies
$$ dF=pFd, \ \ \ Vd=pdV.$$

Here and henceforward, $\s$ denotes the Frobenius map of
$W$ or $W_n$ for all $n$.

We denote by $W\om_Y$ the object of $\A(Y)$,
therefore, a pro sheaf of CDGA's, given by the system of
$W_n\om_Y$'s.
The above identities are often summarized by
saying that $W\om_Y$ is a module for the Raynaud algebra $R=R_0\oplus R_1$,
which is the $W$-algebra generated by symbols $F,V$ in degree
0 and $d$ in degree 1 subject to the relations above.
It should also  be noted that $F$ is multiplicative and agrees
with the usual Frobenius on $W_n \O_Y$.
The Frobenius of $Y$ itself induces a map $\Phi$ of
$W\om_Y$ which is $p^iF$ in degree $i$.

For the purposes of defining the weight filtration, it
will be useful to recall the various different constructions
of the DRW complex. 

One definition of $W_n\om^i_Y$ is as
$$R^iu_{Y/S_n,*}(\O_{Y/S_n, crys})$$
where $u_{Y/S_n}$ is the map from the crystalline site of
$Y$ w.r.t. the base $S_n$ to the \'etale site of $Y$ and $\O_{Y/S_n, crys}$
is the crystalline structure sheaf. That is, 
 $W_n\om^i_Y$ is the $i$-th crystalline cohomology sheaf.
Note here that the cohomology is taken with respect to the
log structure on $S_n$ induced by $L$. 

To `compute' this sheaf, one chooses an {\em embedding system}
$(Y_., Z_.)$ for $Y/S$, that is, a simplicial log scheme $p:Y_.\ra Y$
which is a proper hypercovering for the \'etale topology together
with a closed embedding $Y_. \hra Z_.$, where $Z_.$ is a simplicial
log scheme smooth over $S$. Let $D_.$ be the divided power (PD) envelop
of $Y_.$ in $Z_.$ and let $\Om_{D_.}$ be the associated De Rham
complex. That is, $\Om_{D_n}=\Om_{Z_n/S}\otimes_{\O_Z} \O_{D_n}$ viewed
as a pro-sheaf on $Y_n$, and these
come together to form a simplicial pro-sheaf on $Y_.$ denoted
$C_{Y/S}=C_{Y/S}(Y_,.Z_.)$ and called the {\em crystalline complex} for this
embedding system. Note that $C_{Y/S}$ is actually
a simplicial pro-sheaf of CDGA's. We will denote
by $p$ the map from $Y_.$ to $Y$ viewed as a structure map for
an augmented simplicial scheme and use the same letter
for the maps from the individual
components of $Y_.$. $Rp_*(C_{Y/S})$ is then a cosimplicial
object in the derived category of $\C (Y)$. If we use the natural notation
$s(Rp_*(C_{Y/S}))$ for the complex associated to this cosimplicial
object, then according to \cite{HK} Prop. 2.20,
$$Ru_{Y/S,*}(\O_{Y/S, crys}) \simeq s(Rp_*(C_{Y/S}))$$
canonically as pro-complexes. 
So
$$W\om^i_Y \simeq \H^i (s(Rp_*(C_{Y/S}))).$$
To see in this description the independence of the
embedding system, one merely notes that
any two embedding systems $(Y_.,Z_.), (Y'_.,Z'_.)$
are dominated by a third $(Y_.'',Z_.'')$, which induces
quasi-isomorphisms
$$C_{Y/S}(Y_.,Z_.)\ra C_{Y/S}(Y_.'',Z_.'')
\leftarrow C_{Y/S}(Y'_.,Z'_.)$$

It is worth noting that the DRW complex is a local object,
and that locally, the embedding system can just be
taken as $Y\hra Z$, a lifting of $Y$ to a smooth log scheme over $S$
(which exists by \cite{Ka} Prop. 3.14).
And then we get the simpler formula
$$W\om^i_Y \simeq \H^i (\Om_Z )$$
These formulas just follow from the
definitions, but a  harder theorem (\cite{HK} 4.19)
says that
$$p^*(W\om_Y) \simeq C_{Y/S}$$
as pro-sheaves of CDGA's on $Y_.$ if
we take an embedding system
which admits a 
map $W(Y_.) \ra Z_.$ (for example, if $Z_.$ admits a Frobenius lift). 

This is \cite{HK} Thm. 4.19 except that reference does not
mention the multiplicative structure. However, the compatibility
of this isomorphism with the multiplicative structure
is contained in the proof. More precisely, the map from
$C_{Y/S}$ to $p^*(W\om_Y)$ is determined by
pulling back differentials via the map $W(Y_.)\ra Z_.$
to  the quotient of $\Om_{W(Y_.)/S}$  by the differential graded
ideal
generated by elements of the form
$d(a^{[i]})-a^{[i-1]}da$, where $a\in \mbox{Ker} (W_n(\O_{Y_.}) \ra \O_{Y_.})$,
and then mapping this quotient algebra to
$p^*(W\om_Y)$ by another algebra map. This composed map is
proved to be a quasi-isomorphism.

So we get 
$$TW (W\om_Y)\simeq TW (C_{Y/S}),$$
 which is a q.i. of CDGA's over $K$.

Henceforward, assume that $\Spec(k)$ is equipped with the log structure
of the punctured point, determined by the map $\N \ra k$ of
monoids that sends 1 to 0.

In the construction of the monodromy operator on the crystalline
cohomology of $Y$, a key role is played by the
exact sequence:
$$0\ra W\om_Y[-1] \ra W\tom_Y \ra W \om_Y \ra 0 \ \ \ \ \ \ (*)$$
It is constructed as follows:

Equip $\Spec W[t]$ with the log structure
that takes $1\ra t$.  Let $W<t>$ be the PD envelop of
the ideal $(t)$ in $W[t]$  equipped with the inverse
image log structure.
We have the maps $$\Spec(k)\hra \Spec (W) \hra \Spec (W<t>)$$
which actually are   exact embeddings of log schemes.
Thus, we can consider the crystalline cohomology of $Y$ w.r.t.
$W$ or w.r.t. $W<t>$. On the other hand,
we have the smooth structure map $\Spec W[t] \ra \Spec (W)_0$ where
we use the last subscript 0 to denote that fact
that $W$ is being considered with the trivial log structure.
Let $p:(Y_.,Z_.)\ra Y$ be an embedding system for $Y$ when we
view it as a log scheme  
over $\Spec W[t]$.

 Thus, $Z_.$ can be viewed as a smooth
simplicial log scheme over $\Spec W[t]$ or over $\Spec (W)_0$.
We will abuse notation slightly and denote by the same $Z_.$
the base change of $Z/W[t]$ to $W<t>$.
Therefore, we get a PD pro-De Rham complex
$\Om_{\tilde{D}}$ for the PD envelop of $Y$ in the smooth log scheme
$Z_. \ra \Spec (W)$ (`the De Rham complex of the total space'),
 and
 the pro-De Rham complex $\Om_D$ for the
PD envelop of $Y_.$ in $Z_. \ra \Spec W <t>$ (`the relative De Rham complex'),
which give rise to an exact sequence of pro-sheaves on $Y_.$:
$$0\ra \Om_D [-1] \ra \Om_{\tilde{D}} \ra \Om_D \ra 0$$
The first map takes the differential $\a$ to $\a\wedge dt/t$.
Taking the tensor product with $W$ preserves the exactness (\cite{HK} Lemma 2.22), giving us
$$0\ra \Om_D [-1]\otimes_{W<t>} W \ra \Om_{\tilde{D}}\otimes_{W<t>}W
 \ra \Om_D \otimes_{W<t>}W\ra 0$$
Also, $\Om_D\otimes_{W<t>}W$ is a crystalline complex for
$Y/S$. Now,
$W\tom^i$ is by definition $R^ip_*(\Om_{\tilde{D}}\otimes_{W<t>} W)$.

We  remark that  $W\tom$ carries operators
$F,V,d, \pi $ satisfying the same relations as and compatible
with those of $W\om$: One just repeats the construction of
\cite{HK} p. 246-249 using the $W$ flatness of  $\Om_{\tilde{D}}\otimes W$.

To elaborate a bit on the local description, assume that we have
 a  single smooth log scheme $Z$ over $W[t]$ such that 
$X:=Z\otimes_{W[t]}W$ is a smooth $W$ lift of $Y$. Equip $W[t]$
with the Frobenius which is the usual Frobenius on $W$
and such that $t\mapsto t^p$. Also assume that $Z$
 admits a Frobenius lift $\phi$ (compatible with the Frobenius on
$W[t]$) and that $Z$ is
a lifting of
$Z\otimes_W k$ which is of Cartier type. In particular,
we are assuming that $Z_0 \hra Z$ is an exact closed immersion
of log schemes.
Such liftings always exist locally. In this case,
$\Om_{Z/W}$ and $\Om_{Z/W[t]}$ both carry  operators
$\bF$ which by definition is $\frac{1}{p^i}\phi^*$ on the
$i$-forms ($p^i$-divisibility of $\phi^*$ on
$i$-forms follows from the exactness of the immersion).  
Also,  the exact sequence above becomes
$$0\ra \Om_{X/W}[-1] \ra \Om_{Z/W}\otimes_{W[t]} W \ra \Om_{X/W} \ra 0$$
Here, we have used the fact that, for example,
$$\Om_{\tD/W}=\Om_{Z/W}\otimes_{W[t]}W<t>.$$

We have  the following important formula:
$$(A) \ \ \ \ d^{-1}(p^nC^{i+1})=
\Sigma_{0\leq k\leq n}p^k\bF^{n-k}(C^i)+\sum_{0\leq k\leq n-1}\bF^k(dC^{i-1})$$
where $C$ is either $\Om_{Z/W}$, 
 or $\Om_{X/W}$. 
 The proof is exactly as in \cite{Il}
Prop. 0.2.3.13, where we substitute the Cartier isomorphism
(4.1.2) and the definition (4.3) of \cite{HK} for Illusie's discussion
starting in 0.2.2.2 and up to Prop.2.3.13. This was pointed out by
Jannsen in commentary (6) and (11) to Lemma (1.4) of \cite{Hyo1}.

Formula (A)  implies that
$$\H^i(\Om_{Z/W}\otimes_{W[t]} W_n )\ra \H^i(\Om_{X/W_n)}$$
is surjective for each $n$:
Let $a \in \Om^i_{X/W}$ be closed mod $p^n$. Then
$$a= \Sigma_{0\leq k\leq n}p^k\bF^{n-k}(c^k)+\sum_{0\leq k\leq n-1}\bF^k(db^{k})$$
for some $c^k$ and $b^k$. Lift these to $C^k$ and $B^k$
in $\Om_{Z/W}$. Then
$$A=\Sigma_{0\leq k\leq n}p^k\bF^{n-k}(C^k)+\sum_{0\leq k\leq n-1}\bF^k(dB^{k})$$ is an element of $\Om_{Z/W}$ which is a cocycle mod $p^n$ and maps
to $a$. Now take the image of $A$ in $\Om_{Z/W}\otimes W$
to get the required surjectivity.
Hence, we get the exactness (*) that we want,
as soon as we know that the sequence is independent of
the embedding system. (It also shows that the transition map between
levels is surjective.)
The proof that $W\tom_Y$ and
this sequence is independent of the
embedding system proceeds exactly as the proof of
the independence of the crystalline cohomology: Any two
embedding systems can be dominated by a third, and
the two  sequences possess  maps to the third
which are isomorphisms. That is, once one knows
that the sequences are exact and we have isomorphisms 
for the maps on either end, the middle map is also an isomorphism.

It has been  pointed out to us by the referee
 that the papers \cite{Hyo1} and \cite{HK}
do not deal with the question of whether the exact sequence
which is defined at each level is compatible with
the projection maps. This issue was remedied carefully
in the papers \cite{Nak}, section 6.
Also, in \cite{Nak}, section 10, the Frobenius compatibility
of the exact sequence is checked.

The monodromy operator in crystalline cohomology is the
coboundary map 
$H^i(W\om) \ra H^i(W\om)$
arising from this exact sequence. 

An alternative construction of the monodromy
operator goes as follows:
Define a CDGA $W\tom [u]$ by adjoining to $W\tom$
the divided powers of a variable $u$ in degree zero. That is,
$W\tom[u]$ is the CDGA generated by $W\tom$
and variables $u^{[i]}$, $0, 1, 2, \ldots$,
satisfying the
relations  $du^{[i]}=(dt/t) u^{[i-1]}$ and  $u^{[0]}=1$.
Thus, $W\tom [u]$ is a  `divided-power' Hirsch extension.

Consider the natural map $r: W\tom [u] \ra W\om$ obtained by
composing the two maps $W\tom [u]\ra W\tom \ra W\om$.

\begin{lem} $r$ is a quasi-isomorphism.
\end{lem}

{\em Proof.} Clearly $r$ is surjective.
We can regard $W\tom [u]$ as a double complex whose squares look
like
$$\begin{array}{ccccc}
\ra & W\tom^j u^{[i]} &\ra &W\tom^{j+1} u^{[i]} & \ra \\
 & \da & & \da &  \\
\ra & W\tom^{j+1} u^{[i-1]} &\ra &W\tom^{j+2} u^{[i-1]} & \ra 
\end{array}$$
and we need to prove that the columns are exact for $i\geq 2$.
But the map

$$W\tom^j u^{[i]} \ra W\tom^j u^{[i-1]}$$
takes $\a u^{[i]}$ to $\a dt/t u^{[i-1]}$
and the exact sequence (*) implies the exactness of the sequence
$$\ra W\tom \ra W\tom \ra  W\tom \ra$$
where each horizontal arrow is given by wedging with
$dt/t$, and the columns of our double complex come exactly
from this sequence. This proves the lemma.
\bigskip

$W\tom [u]$ also carries  an $F$-operator by defining
$F(u^{[i]})=p^i u^{[i]}$.

Now define the monodromy operator on $W\tom [u]$ to be
the $W\tom$-linear map that takes $u^{[i]}$ to $u^{[i-1]}$.
The quasi-isomorphism above allows us to transfer this to
$W\om$ in the derived category which hence gives us
a monodromy operator in cohomology. These identities and
$F(dt/t)=dt/t$ give us the relation
$$p\Phi N= N \Phi $$
where $\Phi$ is the operator on $W\om[u]$ and
$W\tom [u]$ equal to $p^iF$ in degree $i$.

\begin{lem} Viewed as an endomorphism of
$R\G(W\om)$ in the
derived category of pro-complexes of $W$-modules,
this monodromy operator agrees with that given
by the coboundary map of the exact sequence (*).
\end{lem}

{\em Proof.}
To compute $N$ on a  cocycle $a$ of $R\G (W\om)$, one lifts it
to a cocycle  $a_0+a_1 u^{[1]}+\cdots \in R\G(W\tom) [u]$, that is,
en element satisfying $[a_0]=a$ and $d(a_0+a_1 u^{[1]}+\cdots)=0$,
applies $N$ to get
$a_1+a_2u^{[1]}+\cdots $ and projects back to $W\om$
to get $[a_1]$. The closedness gives us
$da_0=a_1dt/t$. Thus, $[a_1]$ is the second component of a cocycle 
$(a_0,[a_1])$
in the cone of the map $R\G(W\om [-1]) \ra R\G(W\tom)$ from (*) which
maps to $[a_1]$ via the projection of (*). This proves the claim.
\bigskip

To get the monodromy operator  on the homotopy groups,
we need a slight modification of $W\tom$ to  accomodate augmentations.
For this, it is important to assume that
 $y\in Y$ is a  point  such 
that the log structure on $Y$ is locally of the
form $f^*L$ near $y$, where $L$ is the log structure (of the punctured point)
on $k$ and
$f:Y \ra \Spec(k)$ is the structure map. That is, we assume that
$y:S_0 \ra Y$ is an exact embedding of log schemes.
 
Let $I$ be the kernel of
the natural augmentation map of sheaves
$W\om \ra y_* W.$
$I$ can be resolved by an ideal in $W\tom [u]$
as follows:
Let  $Z$ be a local lifting (around $y$) to a smooth log scheme
over $W[t]$ as above and let $X=Z\otimes_{W[t]} W$. Lift
$y$ to a point $x$ of $X$.
As above, we have the exact sequence
$$0\ra \Om_D\otimes W [-1] \ra \Om_{\tilde{D}}\otimes W \ra \Om_D\otimes W
 \ra 0$$
and we have $\Om_D\otimes W\simeq \Om_X .$ 

There is a natural DG ideal
$\tJ \subset \Om_{\tilde{D}}\otimes W$
defined to be the kernel of the
evaluation at $x$ in degree 0, equal to
$ \Om^i_{\tilde{D}}\otimes W$
for $i\geq 2$, and in degree 1 consists of
the `one-forms whose coefficient of $dt/t$ vanishes
at $x$'.
Only the degree one part requires explanation:
Since  the log structure on $Z$ is locally of the form
$g^*(N)$ near $x$, where $g:Z \ra \Spec(W[t])$ is the structure
map and $N$ is the log structure defined by the divisor $t=0$,
 an element of
$\Om^1_{\tilde{D}}\otimes W$
can be locally written as
$$adt/t +b$$
near $x$ where $b$ is a `usual' regular differential form near $x$.
Therefore, $a(x)$ is well-defined, as is the subsheaf of
1-forms for which $a(x)$ is zero.

Also, denote by $J$ the augmentation
ideal in $\Om_X$ corresponding to evaluation at $x$. Then
an easy local calculation gives us an exact sequence
$$0\ra J[-1] \ra \tilde{J} \ra J\ra 0$$
The cohomology sheaf  of $J$ (regarded as a pro-sheaf on $Y$)
is naturally identified with
the augmentation ideal $I\subset W\om_Y$. Denote by
$\tI$ the cohomology sheaf of $\tJ$.

\begin{lem}
We have a map of exact sequences
$$\begin{array}{ccccccc}
0 \ra & I[-1] & \ra & \tI & \ra & I & \ra 0\\
& \da & & \da & & \da &  \\
0\ra & W\om_Y [-1]& \ra & W\tom_Y & \ra & W\om_Y  & \ra 0 
\end{array}$$
\end{lem}
{\em Proof.}

 We need only check that
the sequence 
$$0 \ra  I[-1]  \ra  \tI  \ra  I  \ra 0$$
is exact. But this follows by the same argument as for the complex
without augmentations.

Now, we can form the complex $\tI[u]$ by adjoining
divided powers of $\log (t)$ as above which gives us
a resolution of $I$. $\tI[u]$ then carries an
operator $N$ again defined by differentiating w.r.t. $u$.

\section{The crystalline fundamental group}

Let $f:Y \ra S_0$ be as in the previous section,
except we further assume that $Y$ is connected and proper. Let
$y \in Y$ be a $k$-rational point which is an exact embedding of
log schemes.

Let
$A_Y:=TW(W\om_Y)$, which we will
call the {\em crystalline rational homotopy type}
of $Y$. This is a CDGA over $K$. It is equipped with
a Frobenius $\Phi$ and a monodromy operator $N$ as described in
the previous section. The point $y$ gives rise to
an augmentation map $a_y:A_Y \ra K$ as follows:
There is a map
$\G (G(W\om_Y)) \ra G(W)$
induced by the map that evaluates elements of $W\O_Y$ at $y$
and sends higher degree elements to zero. (Here, the Godement
resolution $G(W)$ can be taken to be that on the
Zariski site of $W$.) Note that
this is the map induced by the
map of sheaves
$$W\om \ra y_*(W)$$

This induces
a map $A_Y \ra A_k$. $A_k$ consists of compatible collections
of forms on the algebraic $K$ simplices $\D_n$, or in other
words $\invlim \G(\Om_{\D_n})$. Evaluating
 at the zero-dimensional simplex
induces a q.i. $A_k \ra K$, which gives us our augmentation.
From the definition, we see that  the
augmentation ideal is q.i. to
$TW(I)$, where $I$ denotes the kernel of
the map $W\om \ra y_*(W)$.

$B(Y,y):=B(A(Y),a_y)$, is a
 DGA   over $K$ while $H (B (Y,y)))$ is a $K$-CDGA. Now define the $K$-algebra
$Cr(Y,y)$ by the formula
$$Cr(Y,y):= H^0(B(Y,y))$$
Thus, $Cr(Y,y)$ is a $K$-Hopf-algebra and has
a natural filtration $Cr_r(Y,y)=H^0_r(B(Y,y))$, $r=0,1,2,\ldots$
by finitely generated   subalgebras. (Recall here that $H^0_r(B(Y,y))$
is the subalgebra {\em generated} by the $r$-level of the bar filtration.
This is also clearly closed under the comultiplication.)

\begin{definition}
The crystalline fundamental group $\pcr (Y,y)$ is the
pro-algebraic group
defined by  the inverse system
$$\pcr (Y,y)_r:=\Spec (Cr_r(Y,y))$$
\end{definition}
 
\medskip

The $K$-points of $\pcr (Y,y)_r$ are therefore
the group-like elements in the dual Hopf algebra
$Cr_r(Y,y)^*$.

Now assume that $S_0$ has the log structure $L$ of
the punctured point and that the log structure $N$ on $Y$ satisfies
$N_y=f^*(L)_y$.
Recall that we have a quasi-isomorphism from
$TW(\tI [u])$ to $TW(I)$. Thus, we can also compute the
crystalline $\pi_1$ using the bar complex on the
first algebra, which we will denote by
$\TB (Y,y)$. Since the monodromy operator
clearly induces a derivation on $TW (\tI[u])$,
we see that it extends naturally to the tensor product
 derivation on
$\TB (Y,y)$ which commutes with the differential, 
and hence, on $H^0(\TB (Y,y)) \simeq Cr (Y,y)$.
Thus, $N$ is naturally realized as a {\em vector field}
on $\pcr (Y,y)$. Directly from the formulas for the comultiplication
and the fact that $N$ is extended to the bar complex via
tensor products,
we see that $N$ is compatible with the Hopf algebra structure,
and hence, induces a $K$-linear map
 of $\Lie (\pcr (Y,y))$,
since this last is just the dual to the indecomposables
$$QCrys(Y,y):=H^0(B(Y,y))^+/[H^0(B(Y,y))^{+}H^0(B(Y,y))^{+}]$$
Here, the superscript $(+)$ refers to the positively graded part.
Again because it is defined by a tensor product,
$N$ is also compatible with the bar filtration,
hence induces vector fields on all of the
$\pcr(Y,y)_r$.

On $Cr(Y,y)$ as well, the Frobenius and monodromy operators
satisfy the relation
$$p\Phi N=N \Phi$$
since this holds at the level of
the CDGA $W\tom [u]$.

By Proposition 1 of Section 3, we get the relation between
the crystalline fundamental group and crystalline
cohomology:
There is a natural isomorphism
$$H^{cr}_1(Y,K)\simeq \pcr(Y,y)^{ab}$$
where the left hand side denotes the dual
of $H^1_{cr}(Y,K)$, the degree 1 crystalline cohomology
with coefficients in $K$, and the right hand side is the
abelianization of the crystalline fundamental group.
Since the isomorphism is induced by the inclusion
$$I[1] \hra B(Y,y)$$
of the augmentation ideal of $TW(W\om_Y)$ into the
bar complex, it respects the actions of $\phi$
and $N$.

One can also phrase this relation in terms of the
crystalline Lie algebra as
$$H_1(\Lie (\pcr (Y,y))) \simeq H^{cr}_1(Y)^*$$
since the first homology of a Lie algebra
is its abelianization.

\section{The weight filtration in the semi-stable case: proof of
Theorem 1}

We give a few definitions
preliminary to our discussion of the weight filtration.
For this discussion,
let $k$ now be a finite field with $q=p^d$ elements
and let $M$ be an $F$-isocrystal over $k$, i.e.,
a  vector space over $K$ equipped with a $\s$-linear
bijective map $F:M \ra M$.
We say $M$ is pure of weight $i$ if $M$ is finite-dimensional,
 and
the $K$-linear map
$F^d$ has an integral characteristic polynomial
whose roots all have absolute value $q^{i/2}$.
We say $M$ is mixed \cite{Fa} if it has an
increasing (weight-)filtration whose associated graded objects are
pure. We will denote the weight filtration by $P$ (`poids')
because the letter $W$ is used for the DRW complex.

A {\em mixed Frobenius complex } is a triple $(M, N,P)$ of
complexes of $F-$isocrystals $M$ and $N$, where $N$ is
equipped with
an increasing filtration $P$ which is degree-wise
exhaustive and separated
($P_iN_d =N_d$ for $i>>0$ and $P_iN_d=0$ for  $i<<0$)
and a quasi-isomorphism
$M \simeq N$ of $F-$isocrystals, such
that for the spectral sequence
associated to the filtration $P$,
$$E_1^{p,q}=H^{q-p}(Gr_{p} N)$$ 
is a pure $F-$isocrystal of weight $q$. In particular, the
spectral sequence degenerates at $E_2$.
Morphisms  between any of the previous objects
 are required
to respect all structures: that is, the Frobenius, the filtration
and the map between the two complexes $M$ and $N$.

The proof of the following is an easy consequence
of the definitions and the K\"unneth formula.:

\begin{lem} Let $(M_1, N_1, P_1)$ and $(M_2,N_2,P_2)$
be two mixed Frobenius complexes. Then
$$(M_1\otimes M_2, N_1\otimes N_2, P_1\otimes P_2))$$
is a mixed Frobenius complex.
\end{lem}

Define a mixed Frobenius algebra to be a mixed Frobenius
complex $(M, N,P) $ where both $M$ and $N$ are DGA's with
multiplicative Frobenii,
the q.i. $M \simeq N$ respects the algebra structure,
and the filtration $P$ is multiplicative. We will also assume
that all DGA's have connected cohomology.

An augmentation on $(M,N,P)$ is a map of
mixed Frobenius algebras to $(K,K, t)$, where $K$
is the pure crystal of weight zero with $\s$ as Frobenius, and
$t$ is the trivial filtration such that $t_0(K)=K$ and $t_{-1}(K)=0$.
If $(I,J,P)$ denotes the pair of kernels of an augmentation
with the induced filtration, it is clear that $I$ is still quasi-isomorphic
to $J$. Also, since the filtration on $J$ is induced,
we clearly have an inclusion of $E_0$-terms $E^P_0(J)\subset E^P_0(N)$.
Now, suppose we have $[j]\in Gr_k(J)$ such that $[j]=[dn]$ for some
$n \in N$. Since $d$ is $K$-linear, we can always subtract from
$n$ its image under the augmentation. Hence we get
$[j]=[dn']$ for $n' \in J$.

Thus, we conclude that the $E_1$ of $J$ is
a sub-$F-$isocrystal of the $E_1$ for $N$. Therefore,
 $(I,J,P)$ is a mixed Frobenius complex.

Given an augmentation on the mixed Frobenius algebra
$(M,N,P)$, we can construct the bar complexes
$(B(M), B(N), B(P))$, where the filtration $B(P)$ is
the convolution (\cite{Zu}, A.2) of the filtration $P$
and the increasing
bar filtration, that is, $B_n:=\B^{-n}$. 

\begin{lem} If $(M,N,P)$ is a mixed Frobenius algebra,
then $(B(M) , B(N), B(P))$ is a mixed Frobenius complex.
If $(M,N,P)$ is furthermore commutative, then
the bar complex is a mixed Frobenius algebra with
a comultiplication which is a morphism of mixed-Frobenius
algebras. (From the given algebra to the tensor product.)
\end{lem}

{\em Proof.} Let $I$ be the augmentation ideal for $N$.
We examine the spectral sequence for the filtration on
$B(N)$. By the convolution formula (loc. cit.), we see that

$$Gr_{n} B(N)=\oplus_{s+t=n} Gr^P_t Gr^B_s B(N)$$
where $Gr^B$ refers to the graded objects for the bar filtration.
So 
$$\begin{array}{ccc}E_0^{-n,q}&= &\oplus_{s+t=n} Gr^P_t Gr^B_s B(N)^{q-n}\\
 &=&\oplus_{s+t=n} Gr^P_t B(N)^{ -s, q+s-n} \\
&=& \oplus_{s+t=n} Gr^P_t (\otimes^s I)^{q+s-n}\\
&=&\oplus_{s+t=n} Gr^P_t (\otimes^s I)^{q-t}
\end{array}$$
Also one checks readily that the combinatorial differential
is zero on $E_0$, so one gets
$$E_1^{-n,q}=\oplus_{s+t=n}H^{q-t} (Gr^P_t (\otimes^s I))$$
Since $\otimes^sI$ is a mixed Frobenius complex,
we see that this last object is a pure crystal of weight $q$
as desired.

The second sentence follows directly from the formulas
for the multiplication and the comultiplication on
the bar complex (that is, the fact that they preserve
tensor degrees).

\bigskip

The weight filtration in log crystalline cohomology was
studied by Mokrane in \cite{Mok}. The non-trivial
issue of
compatibility of the weight filtrations with
the projection maps was proved in the paper
\cite{Nak} of Nakkajima.

For the purposes of this section, $S_0=\Spec (k)$ is
equipped with the log structure $L$ of the punctured point
and we assume that 
$(Y,M)$ has a log structure which locally fits
into a Cartesian diagram
$$\begin{array}{ccccc}
(Y,M) & \hra & X &\hra &(Z,N) \\
\da & & \da& &\da  \\
S_0 & \hra &S & \hra & \Spec (W[t], N_0)
\end{array} $$
where $Z$ is a $W[t]$ scheme which is smooth over
$W$ with the property that the divisor $X$ over
$t=0$ is of normal crossing with
special fiber $Y$ and $N$ is the log structure
associated to the divisor $E=X+ H$ where $H$
is a divisor which is relatively of normal crossing
meeting $X$ transversally. That is, in \'etale coordinates,
$$Z=\Spec W [t_1, \ldots, t_n]$$
$$X=\Spec W [t_1, \ldots, t_n]/(t_1\cdots t_a)$$
 $E$ is defined by $t_1 t_2 \cdots t_b=0$ for some $b\geq a$,
and the map $Z\ra \Spec (W[t])$ is given by $$t\ra t_1 t_2 \cdots t_a.$$
Also,
$N_0$ is the log structure
associated with the divisor $t=0$. Assume also that $Z$ has a Frobenius lift
$F$ such that $F(t)=t^p(\mbox{unit})$.
Such a $(Z,N)$
is called an `admissible lifting.' 

We also assume that $Y$ itself is globally the union of smooth components
$Y_1, \ldots, Y_c$ which intersect transversally.

At this point we again remind the reader of our convention:
Various pro-sheaves of $W$-modules
will occur in the following discussion.
We will treat them as though they were ordinary sheaves
unless serious confusion is likely to result. This is
especially important to remember in the discussion of
various cohomology sheaves, since, for example,
if $C=\{C_n\}$ is a pro-complex, so that
each $C_n$ is a $W_n$-module, then $\H(C)$ will mean
the pro-system $\{\H(C_n)\}$.

Choose a local admissible lift as above.
Use the notation
$\tom:=\Om_Z(\log E) \otimes_{\O_Z} \O_X$
and $\om:=\tom/(dt/t)$ so that
$W\tom=\H(\tom)$ and $W\om = \H (\om)$.
$\Om_Z(\log E)$ is equipped with the usual weight filtration
$P'$ given by the number of log poles, so that
$\tom $ receives a weight filtration defined by
$$P'_j\tom= Im ((P'_j\Om_Z(\log E) )\otimes \O_X \ra \tom)$$
Therefore,
$W\tom$ is naturally equipped with a weight filtration:
$$P'_j W \tom:=Im [\H (P'_j \tom) \ra W\tom ]$$
 In fact, Mokrane proves that the map
$$\H (P'_j \Om_Z (\log E)\otimes \O_X) \ra W\tom$$
is injective. In \cite{Nak} section 6, it is shown
that the weight filtration defined in this
way for each level is actually compatible with
the projection maps, thereby showing that the filtration
is well-defined on the
pro-complex.

Let $E^i$ denote the disjoint union of the (i+1)-fold
intersections of the fiber over $\Spec(k)$ of
the components of $E$ 
That is:
$$E^i=\coprod_{|I|=i+1} \cap_{j\in I}(E_j\cap Y )$$

In Mokrane \cite{Mok}, where the horizontal component is not considered,
 we had
$$Gr^{P'}_n W\tom =W\om_{E^n}[-n](-n)$$
Here, $[-n]$ refers to a shift of complex degree
while $(-n)$ is a Tate twist for the Frobenius
action. It is important to note here that the $E^i$ are
proper smooth varieties over $k$ and that the graded pieces
are the usual DRW complexes for smooth varieties.

The Hyodo-Steenbrink complex is the simple complex associated
to the double complex
$$WA^{ij}:=W\tom^{i+j+1}/P'_jW\tom^{i+j+1}, \ \ \ i,j \geq 0$$
where the differential
$d':WA^{ij}\ra WA^{i+1,j}$ is induced by $(-1)^jd$, $d$ being
the differential of $W\tom$, and $d'':WA^{ij}\ra WA^{i,j+1}$
is multiplication by $dt/t$.

The map
$$ W \tom \ra W A\ \ \ \a \mapsto \a dt/t$$
factors to a map
$$W\om \ra WA$$
which is proved to be a q.i. Unfortunately, 
this complex is not a CDGA so we cannot
use it to compute the fundamental group.
We also need to take a little care to incorporate
the contributions of the horizontal component.

An alternative construction in the Hodge-theoretic context
was given  in \cite{Ha2}, and
we will use a De Rham-Witt analogue of that construction.
It is based on the simplicial scheme obtained from
the components of $Y$, as occurs often in mixed-Hodge theory.

Let $Y^i$ denote the disjoint union of the (i+1)-fold
intersections of the components of $Y$
That is:
$$Y^i=\coprod_{|I|=i+1} \cap_{j\in I}Y_j$$

For the local admissible lift, we also get the
components $X_j$ of $X$ and the corresponding i+1 fold
intersections $X^i$. We also use the obvious notations
$Y_I$ and $X_I$ for intersections of components indexed by
$I$.
 
Regard the $\O_{X^i}$ as forming a cosimplicial
 pro-sheaf on $Y$.

\begin{lem}
The associated simple complex
$$\O_{X^0}\ra \O_{X^1}\ra \cdots $$
gives a resolution of $\O_X$.
\end{lem}

{\em Proof.} By localizing in the \'etale topology,
we may assume that
the sequence with $\O_X$ appended at the left end
is given by
$$W[t_1,  \ldots, t_n]/(t_1\cdots t_a) \ra 
\prod_i W [t_1, \ldots, t_n]/(t_i) \ra \prod_{i<j}W [t_1,  \ldots, t_n]/(t_i, t_j) \ra \cdots $$
We have canonical isomorphisms
$$W[t_1, \ldots, t_n]/(t_{i_1}, \ldots, t_{i_k}) \simeq
W[t_1, \ldots,  \hat{t_{i_1}}, \ldots, \hat{t_{i_k}}, \ldots, t_n]$$
where the $\hat{\cdot}$ indicates omission of a variable.
Therefore, the complex
$$W[t_1,  \ldots, t_n] \ra 
\prod_i W [t_1,  \ldots, t_n]/(t_i) \ra \prod_{i<j}W [t_1,  \ldots, t_n]/(t_i, t_j) \ra \cdots $$
becomes isomorphic to the tensor product (over $W$) of the complexes
$$W[t_i]\ra W[t_i]/(t_i)$$
for $i=1, \ldots, a$ and the single term complex $W[t_{a+1}, \ldots, t_n]$.
Therefore, the K\"unneth formula tells us that the
only cohomology for this tensor product complex
occurs at the leftmost end, and is  the image of the map
$(t_1)\otimes (t_2)\otimes \cdots \otimes (t_a)$
in $W[t_1, t_2, \ldots, t_n]$, i.e., the ideal $(t_1\cdots t_a)$.
This finishes the proof of the lemma.
\bigskip

\begin{cor}Let $J_i$ be the ideal defining $X_i$
and $J=J_1\cdots J_a$ be the ideal defining $X$. Then
we have an exact sequence of sheaves in the \'etale topology:
$$J_1J_2\cdots J_a \ra \oplus J_i \ra \oplus_{i<j}(J_i+J_j) \ra \cdots$$
\end{cor}
{\em Proof.} Consider the following diagram:
$$\begin{array}{ccccccc}
0  \ra &J & \ra & \O_Z & \ra & \O_X & \ra 0\\
 &\da & & \da & & \da& \\
0  \ra &\oplus_i J_i & \ra & \oplus_i\O_Z & \ra &\O_{X^0}& \ra 0\\
 & \da & & \da & & \da& \\
0\ra &\oplus_{i<j} (J_i +J_j) & \ra & \oplus_{i<j}\O_Z & \ra & 
\O_{X^2}& \ra 0\\
& \da & & \da & & \da& \\
\end{array}
$$
All the rows are exact. The last column is exact by the
previous lemma. The middle column is exact for combinatorial
reasons: The $i-$th level down  is $\O_Z$ tensored with
$i$ copies of the acyclic complex $\Z \ra \Z$.
Therefore, the first column is exact.
\bigskip

Since $\Omega_Z (\log E)$ is flat over $\O_Z$,
 we see that the prosheaves
$$\tom_{X^i}=\tom_{X} \otimes \O_{X^i}=\Omega_Z (\log E)\otimes_{\O_Z}\O_{X^i}$$
also give a cosimplicial resolution of $\tom$.
Here, we say a map from a sheaf $F$ to the degree zero term
of a cosimplicial sheaf
$C$ is a cosimplicial resolution if it yields an exact
complex
$$F\ra s(C)$$
The lemma is saying that
$$\tom \ra \tom_{X^0} \ra \tom_{X^1} \ra \cdots $$
is exact.

We once more remind the reader that in proofs like that of
 the previous
lemma, one should actually be arguing level-by-level, e.g..,
with the rings $W_m[t_1, t_2, \ldots, t_n]$.

\begin{lem} Denote by $C(W\tom_{X})$ the cosimplical
CDGA on $Y$ associated to the collection
$W\tom_{X^i}:=H(\tom \otimes_{\O_Z} \O_{X^i})$,
and let $s(C(W\tom_{X}))$ be the associated
simple complex. So this is the simple complex associated to
the double complex
$$W\tom_{X^0} \ra W\tom_{X^1} \ra \cdots $$
Then
$$W\tom \ra s(C(W\tom_{X})) $$
is a quasi-isomorphism.
\end{lem}

{\em Proof.} Let $J_i$ again be the locally principal ideal
defining $X_i$ in $Z$ and  $J=\prod J_i$.

First, we note that the (inverse) Cartier isomorphism for
$\Om_Z(\log  E)\otimes_W k $ induces one on  $\tom_{X^i}\otimes_k$
for each $i$. To see this, note
the following diagram obtained from that of the previous
corollary by tensoring with $\Om_Z(\log E)$:
$$\begin{array}{ccccccc}
0  \ra &J\Omega_Z(\log E) & \ra & \Omega_Z(\log E) & \ra & \tom & \ra 0\\
 &\da & & \da & & \da& \\
0  \ra &\oplus_i J_i \Omega_Z(\log E) & \ra & \oplus_i\Omega_Z(\log E) & \ra &\tom_{X^0}& \ra 0\\
 & \da & & \da & & \da& \\
0\ra &\oplus_{i<j} (J_i +J_j)\Omega_Z(\log E) & \ra & \oplus_{i<j}\Omega_Z(\log E) & \ra & 
\tom_{X^1}& \ra 0\\
& \da & & \da & & \da& \\
\end{array}
$$
The horizontal sequences are short exact and remain so
after reducing mod $p$.
So we need only show that the Cartier isomorphism on
$\Om_Z(\log  E)\otimes_W k $ induces one on $J\Om_Z(\log  E)\otimes_W k $ and
$(J_{i_1}+\cdots+J_{i_l})\Om_Z (\log E)\otimes k$
for each collection of indices $i_1,\ldots, i_l$.

One sees this  for the complex $J_i \Omega_Z(\log E)\otimes k$ 
by writing it \'etale
locally as the tensor product of the complexes
$$(t_i) \ra k[t_i] dt_i$$ $$k[t_j]\ra k[t_j] dt_j/t_j$$
for $1\leq j \leq b$, $j\neq i$, and
$$k[t_j] \ra k[t_j] dt_j$$
for $j>b$,  noting that the Cartier isomorphism holds for
each term separately (the only new computation is for
the first complex), and using the K\"unneth formula.
Exactly the same argument implies the Cartier isomorphism for
the cohomology of the complex
$$J_{i_1} \cdots J_{i_m}\Omega_Z(\log E)\otimes k$$
for any collection of indices, and hence, for
$J\Omega_Z(\log E)\otimes k$.
For the complex $(J_i+J_j)\Omega_Z(\log E)\otimes k$
we get the result by putting it into an exact sequence
$$0 \ra J_iJ_j \Omega_Z(\log E)\otimes k \ra (J_i \Omega_Z(\log E)\oplus
J_j \Omega_Z(\log E))\otimes k \ra (J_i+J_j)\Omega_Z(\log E)\otimes k \ra 0$$
and using the isomorphism for the first two terms.
One can extend this argument to an arbitrary sum
$$(J_{i_1}+ \cdots +J_{i_m})\Omega_Z(\log E)\otimes k$$
by induction on the number of ideals. That is, the same argument
as in the previous lemma and corollary gives us
an exact sequence
$$\begin{array}{c}
0\ra J_{i_1}\cdots J_{i_m} \ra \oplus_k J_{i_k} \ra \oplus_{k<l} (J_{i_k}+J_{i_l})\\
\ra \cdots \ra J_{i_1}+J_{i_2}+\cdots+J_{i_m}\ra 0
\end{array}$$
and all the terms before the last involve $<m$ ideals. Since the
inverse Cartier map is an  isomorphism  for each of them, it is an
isomorphism for the
last term.

Now we can repeat the definitions of
\cite{Hyo1} p. 245 verbatim to define the maps
${\bf p}: W_n \tom_{X}\ra W_{n+1}\tom_{X}$
and
${\bf p}: W_n \tom_{X^i}\ra W_{n+1}\tom_{X^i}$ for each $i$ and the argument of
\cite{Hyo1} (2.2.2) shows that they are injective.
Furthermore, a repetition of the argument in
\cite{Hyo1} (2.4.1) gives us that
$W_{n+1}\tom_{X}/{\bf p}(W_{n}\tom_{X})$
is q.i. to $\tom_X \otimes_W k$ and
$W_{n+1}\tom_{X^i}/{\bf p}(W_{n}\tom_{X^i})$
is quasi-isomorphic to $\tom_{X^i}\otimes_W k$.
Furthermore, both ${\bf p}$
and the quasi-isomorphisms
$W_{n+1}\tom_{X^i}/{\bf p}(W_{n}\tom_{X^i}) \simeq \tom_{X^i}\otimes_W k$
are compatible with the restrictions
from the complexes on $X^i$ to that on $X^{i+1}$. That is,
we have an exact sequence of double complexes:
$$\begin{array}{ccccccc}
0\ra &W_n \tom_{X}& \ra & W_{n+1}\tom_{X}& \ra & 
W_{n+1}\tom_{X}/{\bf p}(W_{n}\tom_{X}) &\ra 0 \\
 & \da & & \da & & \da & \\ 
0\ra &W_n \tom_{X^0}& \ra & W_{n+1}\tom_{X^0}& \ra & 
W_{n+1}\tom_{X^0}/{\bf p}(W_{n}\tom_{X^0}) &\ra 0 \\
 & \da & & \da & & \da & \\ 

0\ra &W_n \tom_{X^1}& \ra & W_{n+1}\tom_{X^1}& \ra & 
W_{n+1}\tom_{X^1}/{\bf p}(W_{n}\tom_{X^1}) &\ra 0 \\
 & \da & & \da & & \da & \\ 
\end{array}$$
Therefore, the assertion of the lemma is reduced to
the case of $n=1$, in which case it follows again
from the Cartier isomorphism.

\bigskip
We will see below in the proof of independence of the weight filtration
that this resolution also is independent of the admissible lift.
\medskip

A proof identical
to Mokrane's shows that
$$P'_jW\tom_{X^i}:=\H((P'_j(\tom \otimes \O_{X^i}))),$$
where the  $P'$ inside the cohomology is again the
image of the filtration on $\Om_Z(\log E)$,
injects into $W\tom_{X^i}$: 
If one tensors  
the short exact sequence in \cite{Mok} Lemme 1.2 with $\O_{D_i}$
where $D_i$ is one of the components of the
divisor over $t=0$ (which in our notation would
be one of the $X_i$'s), one gets
$$\begin{array}{ccccc}
0\ra P'_{j-1} (\Om_{Z/S} (\log (D))\otimes_{\O_Z} \O_{D_i})
 &\ra& P'_{j}( \Om_{Z/S} (\log (D))\otimes_{\O_Z} \O_{D_i}) & & \\
&\ra &Gr^{P'}_j \Om_{Z/S} (\log (D))\otimes_{\O_Z} \O_{D_i}& \ra& 0
\end{array}$$
Then the retraction used in the proof works equally well
for this sequence to show that the coboundary maps of the
associated cohomology sequence are zero.

For any fixed set $I$ of indices,
the filtration $P'_j$ on $\tom_{X_I}$ is the convolution of
an `internal' weight filtration $P'(I)$ involving log poles coming from
components $X_i$ for $i\in I$,
 and an `external' filtration $P'(E)$
coming from contributions of
$X_k$ for $k\notin I$ and the horizontal divisors
$H_l$'s. Hence, on all of $\tom_{X^i}$ we can break up the filtration
into these two parts.
The graded pieces for the two filtrations 
are calculated separately using the residue formula
as usual:
$$Gr^{P'(I)}_l \tom_{X^i}=\om_{X^{i}}^{\binom{i+1}{l}}[-l]$$
On the other hand, for each indexing set
$I\subset \{1,\ldots, a\}$, we have 
$$Gr^{P'(E)}_k\om_{X_I}=\bigoplus_{J\subset \{1,\ldots, b\}-I, |J|=k }
\Om_{X_{I\cup J}/S}[-k]$$
so we get that
$$Gr^{P'(E)}_k\om_{X^i}\simeq \oplus_K \Om_{X_K}[-k]$$
for various indexing sets $K$ of cardinality $i+1+k$
(allowing multiplicities).
 Therefore, combining  the convolution formula for the
graded pieces  together with the
 the injectivity above,
we see that the terms
$$Gr^{P'}_j W\tom_{X^i}$$
are isomorphic to a direct sum of $W\om_{E_T}[-j](-j)$'s
for indexing sets $T$ of cardinality $i+1+k$ where $k$
runs  between 0 and $j$. Note here that the
$W\om_{E_T}$ are De Rham Witt complexes for smooth
varieties (without log structures).

We can also follow Mokrane essentially literally to
show that the weight filtration $P'$ on each
$W\tom_{X^i}$ is independent of the admissible lifting.
For this reason, we just give a brief sketch:
If $(Z,X, E)$ and $(Z',X', E')$ are two different lifts,
the comparison is effected by constructing the blow-up of
$Z \times_{\Spec(W)} Z'$ along the subscheme 
$\Sigma E_i \times E_i'$ and removing the strict transforms
of each $Z\times E_i'$ and $E_i \times Z'$. Denote the resulting
scheme by $Z''$ and the exceptional divisor by $E''$.
Denote by $X''$ the total transform of just the vertical components
$X_i\times X'_i$.
Finally, let $Y''$ be the intersection with $Z''$
of the total transform of $Y$, embedded diagonally in
$Z\times Z'$.
Thus, the local picture looks as follows:
$$Z=\Spec W[t_1,\ldots, t_n], \ \ \ Z'=\Spec W[t'_1,\ldots, t'_n]$$
$$X=\Spec W[t_1,\ldots, t_n]/(t_1\cdots t_a),\ \ \ X'=
\Spec W[t'_1,\ldots, t'_n]/(t'_1\cdots t'_a)$$
$E$ is given by $t_1\cdots t_b=0$ and $E'$ by
$t'_1\cdots t'_b=0$ for some $b\geq a$.
Then
$$\begin{array}{ccl}
Z''&=&\Spec W[t_1,\ldots, t_n, t'_1,\ldots, t'_n, v_1^{\pm 1}, \ldots, v_b^{\pm 1}]/
(t_1'-v_1t_1, \ldots, t_b'-v_bt_b)\\
&=&\Spec W[t_1,\ldots, t_n, t'_{b+1},\ldots, t'_n, v_1^{\pm 1}, \ldots, v_b^{\pm 1}] \end{array}$$
The exceptional divisor $E''$ is given by
$$t_1\cdots t_b=(\prod v_i)\prod t'_i=0$$
$X''$ by
$$t_1\cdots t_a=(\prod_{i=1}^a v_i)\prod t'_i=0$$
while $Y''$ is given by the ideal
$$(t_i-t'_i , t, p)_{i=1, \ldots, n}$$
Let $\cD'$ be pro-sheaf consisting of
the divided power envelop of $Y''$ in $Z''$
and let $J$ (resp. $J'$) be the smallest
sub-PD ideal of $\cD'$ containing the
ideal defining $Y$ in $Z$ (resp. $Z'$), then
$\cD'/J=\cD'/J'$ and we denote that quotient by $\cD''$.
So $\cD''$ locally looks like
$$\O_{Z''}\otimes_{W[\phi_1, \ldots, \phi_b, \tau_{b+1}, \ldots, \tau_n]}
W<\phi_1, \ldots, \phi_b, \tau_{b+1}, \ldots, \tau_n>$$
where $y_i=v_i-1$ and $\tau_i=t_i-t'_i$.

We wish to compare cosimplicial resolutions of
the three pro-sheaves $\tom_X$, $\tom_X'$ and
$$\tom_{X''}=\Om_{Z''/W} (\log E'') \otimes \O_{X''}$$
If we give the last sheaf as well the filtration $P'$
induced by the filtration on $\Om_{Z''/W} (\log E'')$
by the number of log terms, the above local description
shows that we have natural inclusions
$$P'_j\tom_{X}\hra P'_j\tom_{X''}\otimes_{\O_Z} 
\cD''\hookleftarrow P'_j\tom_{X'}$$
and hence, also inclusions
$$P'_j\tom_{X^i}\hra P'_j\tom_{(X'')^i}\otimes \cD''\hookleftarrow P'_j\tom_{(X')^i}$$
where the superscripts on $X'$ and $X''$ are the obvious ones referring to
tensor products with the structure sheaves of $i$-fold intersections
of the components of $X'$ and $X''$.
But as in  \cite{Mok}, the local description gives
us
$$P'_j\tom_{(X'')^i}\otimes \cD'' \simeq P'_j\tom_{X^i} \otimes_W \Om_{W<\phi_1, \ldots, \phi_b, \tau_{b+1}, \ldots, \tau_n>/W}$$
and 
$$W\simeq \Om_{W<\phi_1, \ldots, \phi_b, \tau_{b+1}, \ldots, \tau_n>/W}$$
so the first inclusion is a q.-i. and similarly for the second inclusion.
If $U$ is the open subscheme of $Y$ obtained by removing the
singular points as well as intersections with the horizontal
divisor $H$,  the cohomology sheaves of all three complexes can
be viewed as subsheaves of the complex $W\tom_U$ from \cite{HK}, section
(1.4). Thus, as explained by Jannsen in the appendix to \cite{Hyo1},
the independence of the filtration follows (that is, no cocycle
condition needs to be checked).

Define a new cosimplicial sheaf 
$C(W\tom [u])$  by adjoining to $C(W\tom [u])$
in each cosimplicial
degree variables
$u^{[i]}$ subject to the condition that $u{[0]}=1$, $du^{[i]}=dt/t u^{[i-1]}$,
each having weight two and Frobenius action $F(u_i)=p^iu_i$.
Thus, $C(W\tom [u])$ gives a cosimplicial resolution of
$W\tom [u]$, that is,
there is a q.i. of complexes
$$W\tom [u] \simeq s(C(W\tom [u]))$$
and hence, a q.i. of CDGA's
$$TW (C(W\tom [u]))\simeq TW(W\tom[u] )$$
On the other hand, we have a natural q.i.
$W\tom [u] \simeq W\om$, and hence, a q.i.
$$TW (C(W\tom[u])) \simeq TW(W\om)$$

We give to $C(W\tom [u])$ and to $s(C(W\tom [u]))$
the weight filtration $P$ defined by the convolution of
the weight filtration on $W\tom [u]$ with
the  filtration by cosimplicial degree: That is,
$$\begin{array}{ccc}
P_n(W\tom_{X^i} [u])&=&P'_{n+i}(W\tom_{X^i}[u])\\
&=&\oplus_j P'_{n+i-2j}(W\tom_{X^i})u^{[j]}
\end{array}$$
By the independence of $P'$ from the lifting, the same
is true for $P$.

As in \cite{Ha2} it is easy to construct a q.i. from
$s(C(W\tom [u]))$ to a cosimplicial resolution of the Hyodo-Steenbrink
complex which is filtration preserving. 
This is achieved by concatenating our q.i.
$W\tom[u] \ra W\om$ with Mokrane's
$W\om \ra WA$ at each level.
Therefore,
on cohomology, we get a filtration preserving isomorphism.
In the case that $k$ is a finite field, we will see
below that the cohomology of $s(C(W\tom [u]))$ is a
mixed isocrystal as was shown for the cohomology of
$WA$ by Mokrane, so a simple strictness argument shows that
the two filtration agree on cohomology.

On the other hand,
$P$ induces a natural filtration on $TW (C(W\tom[u]))$,
and hence, on $B(Y,y)$ and $Cr(Y,y)$.
In general it seems hard to say much about this
filtration. However, the definitions for the situation
where $k$ is finite were concocted to deal with this filtration.

\begin{thm} Suppose the field $k$ is finite.
Then $(TW(W\om), TW (C(W\tom [u])), P)$ is a commutative
mixed Frobenius algebra.
\end{thm}

{\em Proof.}
As filtered complexes, $TW (C(W\tom [u]))$ is filtered q.i. to
$$R\G (s(C(W\tom [u])))$$ so we may compute
the terms of the spectral sequence with the latter.
 We compute the $E_0$ term.
$$\begin{array}{ccc}
E_0^{-n,q}&=&Gr_n [R\G(s(C(W\tom [u])))]^{q-n}\\
&=&\oplus_{0\leq t \leq n} Gr_n [R\G(W\tom_{X^t}[u])]^{q-n-t} \\
&=&\oplus_{0\leq t \leq n } Gr^{P'}_{n+t} [R\G(W\tom_{X^t}[u])]^{q-n-t} \\
&=& \oplus_{0\leq t \leq n} \oplus_i Gr^{P'}_{n+t-2i} [R\G (W\tom_{X^t})]^{q-n-t}(-i)\\
&=& \oplus_{0\leq t \leq n} \oplus_i [G_i[-n-t+2i](-n-t+2i)]^{q-n-t}(-i)\\
&=& \oplus_{0\leq t \leq n} \oplus_i G_i^{q-2n-2t+2i}(-n-t+i)
\end{array}$$
where $G_i$ is a direct sum of 
 complexes $R\G (W\om_{E_T})$ for  some collection of subsets
$T\subset \{1,\ldots b\}$ of cardinality $\geq t$.
Thus,
$$E_1^{-n,q}=\oplus_{t,i}H^{q-2n-2t+2i}(G_i)(-n-t+i)$$
which is pure of weight $q$, being built out of
the crystalline cohomology
of proper smooth varieties. This finishes the proof.

\bigskip
Now given a point $y\in Y$, we get augmentations of
$TW (C(W\tom [u])) $ and $TW(W\om_Y)$. Using them, we can form the
bar complexes
$$B(TW (C(W\tom [u])),a_y)$$ and $B(Y,y)$ which are quasi-isomorphic.
Thus, the weight filtration on $$H(B(TW (C(W\tom [u])),a_y))$$ induces
one on $H(B(Y,y))$ which is compatible with
the Hopf algebra structure.
Therefore, we get a weight filtation on
$Lie\pcr (Y,y)=(QH^0(B(Y,y)))^*$

The previous theorem says that if $k$ is finite, then
$$(TW(W\om_Y), TW (C(W\tom [u])), P)$$ is a mixed Frobenius algebra.
Therefore, we get the structure of a mixed $F-$isocrystal on 
$Lie \pcr (Y,y)$ as well.

\section{Comparison with the de Rham fundamental group: proof of
Theorem 2}

Now suppose $X$ and $Y$ are as in theorem 2.
Denote by $\hX$ the formal completion of $X$ along
$Y$. Regard all of the objects as log (formal) schemes
with the log structure coming from the horizontal divisor
$D$.
Thus, $X$, $\hX$ and $Y$ are all smooth
w.r.t. the trivial log structures on $W$ and $k$.

 We would like
to compare the crystalline fundamental group of
$Y$ with the de Rham fundamental group of
$X^*$, the generic fiber of $X$. At the level of cohomology,
this is the Berthelot-Ogus theorem \cite{BO}
 and we need only to take 
care that the appropriate maps are multiplicative. 
To do this, we need to make explicit several maps in that
paper that are only defined in the derived category by
defining explicit maps between embedding systems that realize
them. 

In the proof, we will be using the computation of the
crystalline fundamental group using the crystalline complex
associated to an embedding system, as explained in section 4.
As a general remark, it is  important to
note that the comparison isomorphism
 between the crystalline complexes associated
to two different embedding systems is effected via a
third embedding  system dominating both, and
the associated inclusion map of the
crystalline complexes is via pull-back of differential
forms. Hence, the crystalline complexes are actually quasi-equivalent
 as algebras. 

We start with  the following preliminary:

\begin{lem} Let 
$$
Y\ra S$$
be a map of fine saturated log schemes over $k$ and let
$S\hra T$ is an exact PD-immersion. Suppose $T$ admits
a Frobenius lift. Then there exists an
embedding system  $(Y_.,Z_.)$
which admits a Frobenius lift compatible with the Frobenius of
$T$.
\end{lem}
{\em  Proof.} Let $\cup_i U_i$ be an affine open covering
of $Y$ such that each $U_i$ lifts to a smooth formal
$T$-log scheme $\tilde{U}_i$ (these exist by \cite{Ka} Prop.3.14).
Let $U=\coprod U_i$, $\tilde{U}=\coprod
\tU_i$. Then the Frobenius $f$ from $U$ to $U$
lifts to $\tU \ra \tU$ by smoothness, in a manner compatible with
the Frobenius of $T$.
On the other hand, $\tU \ra Y$ is  a morphism of
cohomological 2-descent for sheaves having any fixed n-torsion
in the \'etale topology \cite{SD} Cor. 4.3.5. This means
that for the sheaves we are considering (pro-systems   over $n$ of
of $p^n$-torsion sheaves) 
this covering extends (via taking coskeleta)
to a simplicial hypercovering $Y_. \ra Y$ 
such that $Y_i=U\times_Y U \times_Y \cdots \times_Y U$ (i+1 times).
On the other hand, each $Y_i$ embeds naturally  into
an $i+1$-fold  product of $U$  over $S$ which, in turn, embeds
into $Z_i:=\tU\times_T \tU \times_T \cdots \times_T \tU$.
The Frobenius lifting of $Z_1=\tU$ induces
  product maps $Z_i \ra Z_i$ which give us the
required
lifting.
\medskip

Choose a uniformizer $\pi$ of $A$ which therefore
determines a presentation $A\simeq W[t]/(f(t))$,
where $f(t)$ is an Eisentein polynomial of degree $e=[F:K]$.
Let $R$ be the $p$-adic completion of the
divided power envelop of $(f(t),p)$ inside $W[t]$. 
Thus $R$ is also the DP envelop
of the ideal $(t^e)$.
We have a natural map $g:R\ra W<<t>>$. On the other hand,
if $f$ denotes the Frobenius map on $W<<t>>$
which is the usual Frobenius on $W$ and the $p-$th
power map on $t$ and $r$ is such that $p^r\geq e$,
 then the map
$\phi^r:W<<t>>\ra W<<t>>$  factors through
$W<<t>> \ra R \ra W<<t>>$. Note that
the Frobenius map defined above also induces
a Frobenius map $R\ra R$.
Let $Y'=X\otimes A/p$ with the induced log structure.
The Berthelot-Ogus isomorphism hinges on the
comparison between the crystalline cohomology
of $Y$ w.r.t. $W<<t>>$ and $Y'$ w.r.t. $R$.
That is, we have the commutative diagram
$$\begin{array}{ccccc}
Y'& \ra &Y & \hra& Y' \\
\da & & \da & & \da \\
\Spec(A/p) & \ra & \Spec(k) & \hra & \Spec(A/p)\\
\da &  &\da & & \da \\
\Spec(R) & \ra & \Spec (W<<t>>)& \ra & \Spec(R)
\end{array}$$
where the composite of the horizontal arrows are all
the $r$-th iterate of the Frobenius.
We  choose crystalline complexes for $Y$ and $Y'$
as
follows.
Construct first embedding systems for $Y$ and $Y'$
w.r.t. $W[t]$ which fit into a diagram
$$\begin{array}{ccccc}
Y_. & \hra & Y'_. & \hra &Z_. \\
\da & & \da & &  \\
Y&\hra & Y' &  & \\
\end{array} $$
where $Z_.$ is smooth over $W[t]$  and the
left hand square is cartesian. We can also arrange for
$Z_.$ to admit a Frobenius lift compatible with the
Frobenius of $W[t]$. Then
$C=\Om_{Z/W[t]}\otimes W<<t>>$ and
$C'=\Om_{Z/W[t]}\otimes R$ are crystalline complexes
for $Y$ and $Y'$ and we can regard both
as simplicial sheaves of CDGA's on $Y_.$.

Denote by $C^{(r)}$ (resp. $(C')^{(r)}$ ) the pull-back of $C$
(resp. $C'$)
by the $r$-th power of the Frobenius map of $W<<t>>$
(resp. $R$).
Then the big diagram above implies that
$$(C')^{(r)}\simeq C^{(r)}\otimes_{W<<t>>} R$$
in the homotopy category of
sheaves of CDGA's on $Y_.$. On the other hand,
the Frobenius lifts induce maps
$C^{(r)} \ra C$ and $(C')^{(r)} \ra C'$.
So we have  maps  of sheaves of
CDGA's
$$C\otimes R \leftarrow C^{(r)}\otimes R \simeq (C')^{(r)} \ra C'$$
and taking Thom-Whitney algebras, we have multiplicative
maps
$$TW(C)\otimes (R\otimes \Q)  \leftarrow TW(C^{(r)})\otimes (R\otimes \Q)
\simeq TW((C')^{(r)}) \ra TW(C')$$
of CDGA's over $R\otimes \Q$.
The theorem of Berthelot and Ogus imply that
all the maps are quasi-equivalences (since the
usual maps involving the simple complexes rather
than TW algebras are quasi-isomorphisms).
Now we tensor with the quotient map
$R\otimes \Q \ra A\otimes \Q=F$ to get
$$TW(C) \otimes_{W<<t>>\otimes \Q} F \simeq TW(C')\otimes_{R\otimes Q} F$$
By using the fact that
$C'\otimes A$ is the crystalline complex
associated to an embedding system for $Y'$
w.r.t. $A$ which is also true of
$\Om_{\hX/A}$, we get
$$TW(C')\otimes F \simeq TW (\Om_{\hX/A}) \simeq TW (\Om_{X^*/F})$$
where the last quasi-equivalence follows from formal
GAGA.
On the other hand, $C$ is quasi-equivalent to
the base change to $W<<t>>$ of a crystalline complex
for $Y$ w.r.t. $W$ so we get
$$TW(C)\otimes F \simeq TW(W\om_Y)\otimes_K F$$
giving us the desired quasi-equivalence
$$TW(W\om_Y)\otimes_K F \simeq TW (\Om_{X*/F})$$
This gives the isomorphism of homotopy types stated
in the theorem.

It is straightforward to check that this equivalence
is compatible with base-points. So putting everything together,
we see that we get a quasi-equivalence of augmented algebras
$$(TW (\Om_{X^*/K}),x) \simeq  
(TW (W\om_{Y/(W, W(L))},y)\otimes_K F$$
and therefore,
$$Cr(Y,y) \otimes_K F \simeq DR(X^*,x)$$
of the theorem, an isomorphism of commutative Hopf algebras
 over $K$.

As proved in \cite{Wo}, this last object is
the coordinate ring of the De Rham fundamental group of
$X^*$ and the isomorphism is clearly compatible with the
bar filtration (since it is induced by
a q.i. at the level of augmented CDGA's), so we conclude that
$$\pcr (Y,y)\otimes_K F \simeq \pdr (X^*,x)$$

Although we are concentrating on the fundamental group
for this paper, the higher cohomology of the
bar complex is also  of interest.
In particular, if the varieties are simply connected,
they can be used to define higher crystalline
rational homotopy groups.

\bigskip

{\em Proof of Corollary 2.}

The isomorphism classes of the higher rational homotopy groups are
determined by their dimension, and this dimension can be
computed in any complex embedding of $F$ or after base
change to the completion $F_v$ of $F$ w.r.t. $v$. The assumptions
imply that the special fibers $Y$ and $Y'$
 are isomorphic smooth log schemes.
Thus, $TW (\om_Y) \simeq TW(\om_{Y'})$, which implies the
quasi-equivalence of $TW(\Om_X (\log D))\otimes F_v$
and $TW(\Om_X'(\log D'))\otimes F_v$. 
Thus, their bar complexes 
are quasi-equivalent, giving isomorphisms of their cohomology groups,
i.e., the higher De Rham homotopy groups of $U$ and $U'$ \cite{Wo}.
\bigskip

In the projective case without divisors
this theorem follows from Artin and Mazur's
\'etale homotopy theory \cite{AM} where a stronger
integral statement is proved.
The rational statement in the smooth proper case without divisors
 can also be deduced from the formality theorem
of Deligne-Griffiths-Morgan-Sullivan \cite{DGMS} together with
the proper base change theorem for \'etale cohomology and the
comparison theorem between \'etale and Betti cohomologies over ${\bf C}$.

\medskip
{\bf Acknowlegements:} M.K. is grateful to Arthur Ogus
for several useful discussions related to
the comparison theorem and to
Kirti Joshi for numerous conversations on p-adic Hodge theory.
He is also grateful to the Korea Institute for Advanced Study
for providing the environment
where the revisions to this paper could be carried out.

The referee of an earlier version of this paper gave it
an unusually thorough reading and made many helpful
suggestions for improving the proofs. For
this we are extremely grateful.

Both authors were supported in part by grants from the
National Science Foundation.

{\footnotesize
M.K.: DEPARTMENT OF MATHEMATICS, UNIVERSITY OF ARIZONA, TUCSON, AZ 85721, EMAIL: kim@math.arizona.edu and
KOREA INSTITUTE FOR ADVANCED STUDY, 207-43 CHEONGRYANGRI-DONG, DONGDAEMUN-GU,
SEOUL, KOREA}

{\footnotesize R.H.: DEPARTMENT OF MATHEMATICS, DUKE UNIVERSITY, DURHAM, NC 27708, EMAIL: hain@math.duke.edu}

\end{document}